\documentclass [a4paper,12pt] {amsart}
\usepackage [latin1]{inputenc}
\usepackage[all]{xy}
\usepackage{amsmath}
\usepackage{amsfonts}
\usepackage{amssymb}
\usepackage{graphics}
\usepackage[lite]{amsrefs}
\usepackage{amsfonts}
\usepackage{amssymb}
\usepackage{xy,amsmath,amssymb}
\usepackage[dvips]{graphicx}
\usepackage{graphicx}
\usepackage[lite]{amsrefs}
\newtheorem{theorem}{Theorem}[section]
\newtheorem{lemma}[theorem]{Lemma}
\newtheorem{prop}[theorem]{Proposition}
\newtheorem{cor}[theorem]{Corollary}
\newtheorem{ex}{example}[section]
\newtheorem{remark}{Remark}[section]
\newtheorem{defin}[theorem]{Definition}
\newtheorem{digr}[theorem]{Digression}
\def\be{\begin{equation}}
\def\ee{\end{equation}}
\def\bt{\begin{theorem}}
\def\et{\end{theorem}}
\def\bp{\begin{prop}}
\def\ep{\end{prop}}
\def\bl{\begin{lemma}}
\def\el{\end{lemma}}
\def\bc{\begin{cor}}
\def\ec{\end{cor}}
\def\br{\begin{oss}\rm}
\def\er{\end{oss}}
\def\bex{\begin{ex}\rm}
\def\eex{\end{ex}}
\def\bd{\begin{defin}}
\def\ed{\end{defin}}
\def\demo{\par\noindent{\bf Proof.\,\,}}
\def\bit{\begin{itemize}} \def\eit{\end{itemize}} 
\def\enddemo{\ $\Box$\par\vskip.6truecm}
\def\R{{\mathbb R}}   \def\a {\alpha} \def\b {\beta}\def\g{\gamma}
\def\N{{\mathbb N}}     \def\d {\partial} 
\def\C{{\mathbb C}}
\def\D{{\mathbb D}}     
                 \def\p{\partial}\def\r{\varrho}
        \def\z{\zeta}
 \def\tms{\times} 
\def\sbs{\subset} \def\Sbs{\Subset} 

 \def\oli{\overline} 

 \def\uli{\underline} \def\oli{\overline}

  \def\beq{\begin{eqnarray}}
\def\eeq{\end{eqnarray}}

\newcommand{\ocal}{\mathcal{O}}
\newcommand{\hcal}{\mathcal{H}}
\newcommand{\ycal}{\mathcal{Y}}
\newcommand{\xcal}{\mathcal{X}}
\newcommand{\zcal}{\mathcal{Z}}
\newcommand{\fcal}{\mathcal{F}}
\newcommand{\gcal}{\mathcal{G}}
\newcommand{\ucal}{\mathcal{U}}
\newcommand{\vcal}{\mathcal{V}}
\newcommand{\wcal}{\mathcal{W}}
\newcommand{\ccal}{\mathcal{C}}
\newcommand{\dcal}{\mathcal{D}}

\newcommand{\bcal}{\mathcal{B}}
\newcommand{\tcal}{\mathcal{T}}
\newcommand{\scal}{\mathcal{S}}
\newcommand{\pcal}{\mathcal{P}}
\newcommand{\ical}{\mathcal{I}}

\newcommand{\homs}{{\sf Hom}}

\newcommand{\hsolo}{\mathcal{H}_s}

\newcommand{\hsolop}{{\mathcal{H}_s}_\bullet}

\newcommand{\ocolim}{\mathrm{hocolim}}
\newcommand{\olim}{\mathrm{holim}}

\newcommand{\xdot}{[X_\cdot]}
\newcommand{\zdot}{[Z_\cdot]}
\newcommand{\pdot}{[P_\cdot]}

\newcommand{\drm}{\mathrm{d}}

\newcommand{\fdot}{F_\cdot}

\newcommand{\pis}{\pi^{\text{simpl}}}

\newcommand{\cbb}{\mathbb{C}}
\newcommand{\dbb}{\mathbb{D}}
\newcommand{\rbb}{\mathbb{R}}
\newcommand{\nbb}{\mathbb{N}}

\newcommand{\fsf}{{\sf F}}
\newcommand{\gsf}{{\sf G}}
\newcommand{\csf}{{\sf C}}

\newcommand{\pref}{{{\sf Prsh}}_T(\mathcal S)}

\newcommand{\etale}{\'etale\;}

\newcommand{\cce}{\'e\;}
\newcommand{\acc}{\`a\;}

\def\ssmi{\!\smallsetminus\!}

\title[ANALYTIC STACKS AND HYPERBOLICITY]{ANALYTIC STACKS AND HYPERBOLICITY}
\author[Simone Borghesi and Giuseppe Tomassini]{Simone Borghesi and Giuseppe Tomassini}
              \address{Universit\acc degli Studi di Milano - Bicocca -  Via Cozzi, 53 - 20125 Milano, Italy}
              \email{simone.borghesi@unimib.it}
              \address{Scuola Normale Superiore, Piazza dei Cavalieri, 7 - I-56126 Pisa, Italy}
              \email{g.tomassini@sns.it}
 \date{\today}
            
\begin {document}
 \maketitle
\hfill{\small{\it Dedicated to Franco Bassani}}          
\vskip45pt 
\begin{abstract}
In this manuscript we give two notions of hyperbolicity for groupoids on the analytic site of complex spaces; 
we call them Kobayashi and Brody's. They generalize the classical ones, in the particular case the groupoid 
is a complex space. We prove that such notions are equivalent if the groupoid is  a {\em compact} Deligne-Mumford
analytic stack (in analogy with the Brody Theorem). Moreover, under the same assumptions, such notions of 
hyperbolicity  are completely detected  by the coarse moduli space of the stack. We finally show that stack hyperbolicity, as we defined it, is expected to impose a peculiar behaviour to the stack itself, 
much like hyperbolicity for complex spaces; 
for instance, a stronger notion of it (hyperbolicity of the coarse moduli space) implies a ``strong asymmetry'' on the stack, its automorphism 2-group having only finitely many isomorphism classes, in the compact case.
\end{abstract}
 \tableofcontents
\section{Introduction}
A classical notion, introduced by Kobayashi, is the one of a 
pseudodistance associated with any complex space which behaves as a contraction with respect to 
holomorphic maps \cite{KOB1}. In particular, this pseudodistance is a biholomorphic invariant. Hyperbolic complex 
spaces are precisely those for which this pseudodistance is, in fact, a distance.  
In this manuscript, a {\em Kobayashi hyperbolic space} is what we referred to as a hyperbolic complex space, whereas a {\em Brody hyperbolic space} is a complex space not admitting nonconstant holomorphic maps from $\cbb$; the latter notion is inspired by Brody's Theorem \cite{brody}. 
For the basic results in hyperbolicity, implications and conjectures in complex geometry we refer to 
\cite{KOB1}. See also \cite{brun}.

At the core of the paper lies the notion of  hyperbolicity for analytic stacks following this Brody semantics. Such notions ought at least to generalize the known ones for complex spaces and be categorical equivalence invariant (thus invariant by the choice of the presentation of an analytic stack). 

The exposition of the manuscript begins in a greater generality by dealing with $\scal$-groupoids (groupoids for short), where $\mathcal S$ is the category of (Hausdorff and reduced) complex spaces and coverings those induced by the strong topology, and it unfolds by progressively adding supplementary assumptions to the studied objects, when necessary,  and restricting to Deligne-Mumford analytic stacks for the proof of the main theorems.

In the first two sections we list the necessary homotopical results for providing insight to the hyperbolicity definitions.
Some of the notions appearing in the paper \cite{BT2} in terms of simplicial sheaves can be effectively expressed in the category of groupoids by using the results appeared in \cite{DHI} and \cite{HOL}. The simplicial presheaf viewpoint enables 
us to adapt to groupoids the concept of Brody hyperbolicity already introduced in the paper \cite{BT2}, whose motivation we briefly recall here. 
The Brody hyperbolicity condition for a complex space $Y$, summarized in the bijectivity of
$$
p^*:\homs_{\rm holo}(X, Y)\to \homs_{\rm holo}(\cbb\times X, Y)
$$
for all complex spaces $X$ where $p:\cbb\times X\to X$, can be extended to groupoids $\gcal$ in two ways. One by requiring the bijectivity of $$
p^*:\homs_{{\rm Grp}/\scal}(\xcal,\gcal)\to  
\homs_{{\rm Grp}/\scal}(\cbb\times \xcal, \gcal)
$$
for all complex spaces $\xcal$ and the other for all groupoids $\xcal$. In analogy with the classical definition we choose the latter 
(see Definition \ref{brody-iper1}). With some more work, we rephrased this definition in terms of {\em presheaves of holotopy groups} 
of the groupoid (see Definition \ref{br-campi}). In this form it is clear that groupoid Brody hyperbolicity extends the same classical property for complex spaces and that it is categorical equivalence invariant. Contradditing Brody hyperbolicity of an analytic stack provides a way to showing the existence of $\C$ parametrized families of objects in the moduli problem associated with the stack.


In order to introduce Kobayashi hyperbolicity for groupoids $\gcal$, several definitions have been considered; the one we decided to 
use (see Definition \ref{kob-iper}) it involves the holotopy presheaves $\pis_0$, $\pis_1$, as well. It is based on the 
notion of {\em relative analytic disc} and {\em relative analytic chain} (see \ref{dischi}) joining two sections in $\pis_0(\gcal,g)(U)$ or $\pis_1(\gcal,g)(U)$, $U$ being a complex space. It preserves the metric ``flavour'' of the classical notion, bearing a difference: two points of a connected complex space are always joined by an analytic chain, whereas for the holotopy presheaves this happens only for particular sections called {\em admissible}. For a pair of admissible  sections of $\pis_0(\gcal,g)(U)$ or  $\pis_1(\gcal,g)(U)$ a {\em Kobayashi pseudodistance} ${\rm d}_{\rm Kob}$  is defined (see \ref{ddisc}).
If this is strictly positive for any pairs of admissible sections of $\pis_0(\gcal,g)(U)$ and $\pis_1(\gcal,g)(U)$
over any complex space $U$, the groupoid $\gcal$ is said to be {\em Kobayashi hyperbolic}.

The important Brody's Theorem, 
\cite{brody}, states that a compact complex space is hyperbolic if the only holomorphic maps from $\cbb$ to it are 
constant. 
One of the main results of the paper is the generalization of Brody's theorem to compact Deligne-Mumford analytic stacks (compactness for stacks is discussed in Section \ref{aus}).

The proof takes place in the Sections \ref{kob-sec} and \ref{dura} and it extensively uses techniques of complex variable, made possible by the replacement of $\ycal$ by its presentation 
$$[X_1=X\times_\ycal X\rightrightarrows X]$$
since the properties we will be studying are equivalence invariant. 
A crucial role is played by the {\em coarse moduli space $Q(\ycal)$} of a Deligne-Mumford analytic stack $X\to \ycal$, whose points are equivalence classes of an equivalence relation on $X$. The geometrization of this set was proved, in a more general context, in \cite{grezzo} and is a key step to make all the parts of the proof patch together. In our case, by fixing a distance function on $Q(\ycal)$ induced by a length function and lifting it to $X$ and $X_1$, we obtain distances such that the structure maps are local isometries. This allows us to define distances $\delta_U$, ${\delta}_{\pi_0,U}$, ${\delta}_{\pi_1,U}$ on the sets ${\sf Dis}_{\xdot}(U)$, $\pi_0(\ycal,y)(U)$ and $\pi_1(\ycal,y)(U)$ (see Proposition \ref{distp2}) and show the following fundamental estimate (see Lemma \ref{lcat1}): for every complex space $U$ there exists a positive number $c(\ycal;U)\le+\infty$ such that if $\alpha_1,\alpha_2\in \pis_i(\ycal,y)(U)$, $i=0,1$, then
$$
 {\rm d}_{\rm Kob}(\a_1,\a_2)\ge\frac{{\delta}_{\pi_i,U}(\a_1,\a_2)}{2\, c(\ycal;U)}
$$
Thus Kobayashi hyperbolicity of a stack $\ycal$ is reduced to the finiteness of $c(\ycal;U)$ for any complex
space $U$. 

The statement ``Kobayashi hyperbolicity implies Brody hyperbolicity'' holds in general. This fact, which in the classical case is a simple consequence of non Kobayashi hyperbolicity of $\cbb$ and that every holomorphic map is a contraction with respect to the Kobayashi presudodistance, in the context of stacks is not entirely obvious (see Theorem \ref{KB}). 

In Section \ref{dura} we prove that for a Deligne-Mumford analytic stack the converse is also true (see Theorem \ref{BRO}). We argue by contradiction assuming that $c(\ycal;U)=+\infty$ for some complex space $U$; then, by the results of Section \ref{top}, there is a sequence $\{f^\nu\}_\nu$ of holomorphic 
maps $f^\nu:\dbb\to Q(\ycal)$ such that $\lim_{\nu\to +\infty}|\drm f^\nu(0)|=+\infty$. By classical results, such as 
the ``reparametrization Lemma'' (cfr. \cite{brody}) and Ascoli-Arzel\acc Theorem, there exists a subsequence 
$\widetilde f^\mu$ uniformly convergent on compact sets to a holomorphic map $f:\cbb\to Q(\ycal)$, which is not constant since $|{\rm d}f(0)|=1$. The hard part of the proof consist in showing that actually $f$ lifts to a non constant morphism $\C\to\ycal$.

The techniques previously developed  allow to prove  the second main theorem of the paper, namely that for a compact Deligne-Mumford analytic stack $\ycal$, hyperbolicity follows from the classic hyperbolicity of the complex space $Q(\ycal)$ (see Corollary \ref{equivalenze}).

Finally, in Section \ref{hyp-coars} we give some applications of the methods and results contained in the previous sections. We show that hyperbolicity for compact Deligne-Mumford stack is not equivalent to the hyperbolicity
of its coarse moduli space, by providing an explicit example of an hyperbolic Deligne-Mumford stack with a torus as 
moduli space. Furthermore, we prove that the
hyperbolicity of the coarse moduli space implies that the compact Deligne-Mumford stack has few automorphisms; more precisely, the isomorphism classes of the $2$-group ${\sf Aut}(\ycal)$ are finitely many (Theorem \ref{auto}).

We wish to thank Gabriele Vezzosi for several discussions regarding the homotopic part of the paper, Angelo Vistoli who provided the core steps in the proof of Theorem \ref{auto}, Giorgio Ottaviani, Jean Pierre Demailly and 
Burt Totaro whose comments helped to better understand to role of hyperbolicity of the coarse moduli space.


\section{Preliminaries}
 
\subsection{Notation and definitions}\label{notation}
\begin{itemize}
\item $\scal_T$ is the analytic site: the category $\mathcal S$ whose objects are complex spaces and coverings those induced by the strong topology.
\item ${\rm Grp}$ is the category of (set-theoretic) groupoids and by ${\sf Psh}(\mathcal S;{\rm Grp})$ we will denote the category of presheaves on $\scal$ of set-theoretic groupoids. 
\item ${\rm Grp}/{\mathcal S}$ is the category of $\mathcal S$-groupoids whose objects are categories fibered in groupoids. A stack is an object of ${\rm Grp}/{\mathcal S}$ which is a sheaf on $\scal_T$ 
\item $\Delta^{op}\pref_J$ is the category of simplicial presheaves of sets on the site $\scal_T$ with a topology $T$ and endowed by {\em local, injective, simplicial model structure} on $\Delta^{op}\pref$ (cfr. Joyal's model structure, \cite[ Section 5.1]{HOL}).  $\hsolo$ (respectively $\hsolop$) is the homotopy category (respectively the pointed homotopy category) associated.
\item Let $\xcal$ be a groupoid. Then $c_\xcal:\xcal\to \ccal\xcal$ is its campification (cfr. \cite[Lemma 3.2 and Observation 3.2.1  (3)]{lm}).
\item The morphisms $\partial_0$ and $\partial_1$ will denote the face morphisms of a simplicial presheaf $\xcal_\cdot$ or, more frequently, of groupoids $\xdot$ and $\ccal\xdot$ from the presheaf in degree $1$ to presheaf.
\end{itemize}\bigskip


\subsection{Simplicial presheaves and groupoids}

Before being able to state what we think of as a (differently flavoured) hyperbolic groupoid, we will 
expose the connection between simplicial presheaves and groupoids, by stating the relevant results contained 
in the manuscript \cite{DHI}. The concept of Brody hyperbolicity, in particular, is directly transposed from 
simplicial presheaves. While the concept of groupoids in terms of categories fibered in set-theoretic groupoids 
($S$-groupoids) probably goes back to ideas of Grothendieck, only recently these objects have been related  
to the homotopy theory of simplicial presheaves of sets.

In the paper \cite{BT2} we dealt with simplicial presheaves of sets on the analytic site, i.e. the category $\Delta^{op}\pref$ whose 
objects are complex spaces and coverings given by the strong topology on them. That category has been considered a model category by means of the Joyal injective {\em local} simplicial structure.  Here the word ``local'', as opposed to ``global'', refers to the fact that the weak 
equivalences are morphisms inducing simplicial sets weak equivalences on the stalks of the presheaves, as 
opposed to weak equivalences on simplicial sets of sections of the simplicial presheaves. 
After having localized this category with respect of such model structure, we defined Brody hyperbolic 
simplicial presheaves (cfr. \cite[Section 3.1]{BT2}). 
We want to relate these objects to groupoids in general and stacks in particular. 

Such a relation is spread out in the papers \cite{HOL} and \cite{DHI}. We list here their results which are relevant to this manuscript.  
The starting point is Corollary 4.3 in \cite{HOL} which claims the existence of an adjunction between the 
categories ${\sf Prsh}({\rm Grp})$ and ${\rm Grp}/{\mathcal S}$. Endowing each of appropriate {\em global} model structures (\cite[Proposition 4.1 and Theorem 4.2]{HOL}) we furthermore have that the adjoint functors are a Quillen equivalence. Here 
``global'' means that the weak equivalences are meant to be objectwise (respectively fiberwise) weak 
equivalences. These two categories are not directly related to $\pref$, but a certain localization of them are.
Keeping the same notation as the references, we let $S$ to be the set of maps

\begin{equation}
S=\{p_U: \mathrm{hocolim}\, U_\cdot \to X:\; \amalg_iU_i= U\to X\; \text{is a cover for the strong topology}\} 
\end{equation} 
where $\mathrm{hocolim}$ is the homotopy colimit of the diagram $U_\cdot=\{\cdots U\times_X U\rightrightarrows U\}$.
By using Bousfield localization theory, it can be shown that there are model structures on ${\sf Prsh}({\rm Grp})$ and ${\rm Grp}/{\mathcal S}$ whose homotopy category is the localization with respect of $S$ (\cite[Proposition 4.4]{HOL}). 
These model structures are called {\em local} for a reason that will become clear later and to explicit this 
model structure we will add $L$ as a subscript. It follows that the aforementioned adjoint pair of functors, 
induces a Quillen equivalence ${\sf Prsh}({\rm Grp})_L\leftrightarrows ({\rm Grp}/{\mathcal S})_L$ (\cite[Corollary 4.5]{HOL}).
In the paper \cite{DHI} the model structure on $\Delta^{op}\pref$ is reinterpreted through the Bousfield 
localization with respect to the maps in $S$: by \cite[Theorem 6.2]{DHI}, we know that it localized category 
is equivalent to the homotopy category $\hsolo$. Thus, to relate $\Delta^{op}\pref_J$ to $({\rm Grp}/{\mathcal S})_L$, it 
suffices to relate $\Delta^{op}\pref_L$ to ${\sf Prsh}({\rm Grp})_L$. Again, there is an adjunction  $(\pi_{oid}, N)$
$$
\pi_{oid}: \Delta^{op}\pref\leftrightarrows {\sf Prsh}({\rm Grp}):N
$$ 
defined as follows: $\pi_{oid}$ is the functor which sends a simplicial presheaf $\fdot$ to the $\mathcal S$-groupoid
having $F_0(U)$ as objects, for each complex space $U$, and $\homs(a,b)$ is the set of $\phi\in F_1(U)$
such that  $\partial_0(\phi)=a$ and $\partial_1(\phi)=b$.
To an $\mathcal S$-groupoid  $\gcal$, the functor $N$ associates the simplicial presheaf $N\gcal$ with $(N\gcal)_0={\sf Ob}
(\gcal)$, $(N\gcal)_1={\sf Mor}(\gcal)$ and 
$$
(N\gcal)_i=(N\gcal)_1\times_{(N\gcal)_0}\stackrel{i}{\cdots}\times_{(N
\gcal)_0} (N\gcal)_1
$$
with the following structural face morphisms: $\partial_0,\partial_1: (N\gcal)_1\to (N
\gcal)_0$ are the domain and codomain of the isomorphism, respectively;  the three morphisms $(N\gcal)_2\to 
(N\gcal_1)$ send $(f,g)$ respectively in $f$, $g\circ f$ e $g$; in the general case an $n$-tuple of composable 
isomorphisms are sent to $(n-1)$-subtuples involving individual isomorphisms and compositions of them, 
when appropriate. The degenerations are induced by alternatively adding the identity morphism.
$(\pi_{oid},N)$ is a Quillen pair and the following holds (cfr. \cite[Theorem 5.4]{HOL}):

\begin{theorem}\label{null}
The Quillen adjoint pair $(\pi_{iod},N)$ induces a Quillen equivalence between $(\mathbb S^2)^{-1}\Delta^{op}\pref_L$, 
the $\mathbb S^2$ nullification of $\Delta^{op}\pref_L$, and ${\sf Prsh}({\rm Grp})_L$.
\end{theorem}

To trace this chain of equivalences back to the category $\Delta^{op}\pref_J$ equipped with the Joyal simplicial model structure, we finally need to use results in the \cite{DHI} which imply that $\Delta^{op}\pref_J$ and $\Delta^{op}\pref_L$ are Quillen equivalent.
The previous considerations and this theorem prove

\begin{cor}\label{null-import}
The adjoint pair $(\pi_{oid}, N)$ induces a Quillen equivalence between the categories $(\mathbb S^2)^{-1}\Delta^{op}\pref_J$ and $({\rm Grp}/{\mathcal S})_L$.
\end{cor}

In order to see how this relates to (analytic) stacks we need to invoke sharper results. Theorem 1.1 (see also Theorem 3.9) of \cite{HOL} states that a groupoid $\fcal$, seen as a presheaf of groupoids, is a stack if and only if, for any covering $\amalg_iU_i=U\to X$ , the canonical morphism 
$$
\fcal(X)\to \text{holim}_n \Big\{\prod \fcal(U_i)\Rightarrow \prod \fcal(U_{ij})
\Rrightarrow \prod \fcal(U_{ijk})\cdots\Big\}
$$
is an equivalence of categories for each complex space $X$, where $U_{i_1,\cdots,i_k}$ stands for \\$U_{i_1}\times\cdots\times U_{i_k}$. 
On the other hand, since the Bousfield $\mathcal S$-localizing structure on ${\sf Prsh}({\rm Grp})$ relies on an underlying {\em global} (meaning weak equivalences and fibrations are objectwise) model structure, we deduce that the $\mathcal S$-fibrant objects are precisely those presheaves of groupoids $\mathcal F$ that are
\begin{itemize}
\item objectwise fibrant, i.e. all since a set theoretic groupoid is simplicially fibrant if seen as simplcial set (by means of the functor $N$, precedently described) and
\item such that the canonical morphism 
$$
\mathcal F(X)\to \mathrm{holim}_n \Big\{\prod \mathcal F(U_i)\Rightarrow \prod \mathcal F(U_{ij})\Rrightarrow \prod \mathcal F(U_{ijk})\cdots\Big\}
$$
is a weak equivalence of set theoretic groupoids.
\end{itemize}

Because of Theorem \ref{null}, the functor $\pi_{oid}$ sends $\scal$-fibrant presheaves of groupoids to $\mathcal S$-fibrant simplicial presheaves and by \cite[Theorem 1.1]{DHI} such simplicial presheaves are precisely those that are fibrant according to the Joyal simplicial model structure on $\Delta^{op}\pref$. This leads to 

\begin{theorem}\label{fibranti}
The chain of Quillen equivalences between ${\rm Grp}/\mathcal S$ and $(\mathbb S^2)^{-1}\Delta^{op}\pref_J$ induces an 
isomorphism between the (full) subcategories of stacks and fibrant simplicial presheaves. 
\end{theorem}

While for $\scal$-groupoids there are two notions of equivalences, one {\em global} and one {\em local}, they coincide
for stacks. Moreover, we may think of a stack as a simplicial presheaf of sets, presheaf of groupoids, or categories
and equivalences as homotopy equivalences in the first two cases and category equivalences in the latter.

\subsection{Analytic stacks}\label{campi-alg}

\noindent A stack over the analytic site $\scal_T$ is said to be an {\it analytic stack} if\vspace{2mm}
 \bit
 \item[1)] the diagonal morphism $\Delta_{\ycal}:\ycal\to\ycal\times\ycal$ is representable; 
 \item[2)] there exists a complex space $X$ with a smooth and surjective morphism $p:X\to \ycal$.
 \eit
 \vspace{2mm}\noindent
The morphism $p:X\to \ycal$ is also called a {\em presentation} of $\ycal$ and $X$ an {\em atlas}. 

An analytic stack ${\ycal}$ with an \etale presentation $p:X\to \ycal$ is called a {\it Deligne-Mumford analytic stack}.

Let $\pcal$ be a presheaf of grupoids and $F:\pcal\to\gcal$ be a $1$-morphism (functor) to a groupoid $\gcal$. We build a groupoid out of it, denoted by $[\pcal_\cdot]$. Its objects over a complex space $U$ are the sections in $\pcal(U)$ and $\homs_{\pdot(U)}(f,g)$ are the sections $\phi\in \pcal_1(U):=\pcal(U)\times_\ycal \pcal(U)$ such that $\partial_0(\phi)=f$ and $\partial_1(\phi)=g$, where $\partial_i:\pcal_1\to \pcal$, for $i=0,1$, are the projections on the factors, and the fiber product is taken in the category of groupoids. The remaining structure making $\xdot$ a groupoid is inherited by the one of $\gcal$ and it is explained in (2.4.3) and in Proposition 3.8 of \cite{lm}. 

\begin{remark}\hspace{1cm}
{\em
\begin{itemize}
\item[1)] If $X\to \ycal$ is an analytic stack (see Section \ref{notation}), then the objects of $\xdot$ over a complex space $U$ are the holomorphic maps $U\to X$ and $\partial_i:X_1\to X$ are holomorphic maps between complex spaces.  
\item[2)] Our notation is slightly different from the one in \cite[2.4.3]{lm}: the groupoid $\xdot$ is denoted as $\xdot'$ there. Moreover, through this manuscript, we will identify the {\em $\mathcal S$-espace en groupo\"ides} and its associated groupoid.  
\end{itemize}
}
\end{remark}

We recall the following general result
\begin{prop}(cfr. \cite[Prop. 3.8]{lm}) Let $F:\pcal\to\ycal$ be a morphism (functor) between a presheaf and a stack. Then the canonical morphism (functor) $[\pcal_{\cdot}]\to\ycal$ is a monomorphism and is epi if and only if $F$ is.
\end{prop}

In the particular case $p:X\to\ycal$ is an analytic stack, we get a simplicial complex space $\xdot$ such that 
$$
X\times_\ycal\times\stackrel{i+1}{\cdots}\times_\ycal X=X_i=X_1\times_{\partial_0,X,\partial_1}\stackrel{i}{\cdots}\times_{\partial_0,X,\partial_1} X_1.
$$
$\xdot$ is a {\em prestack}, as explained in the example \cite[3.4.3]{lm} and is precisely $N([\pcal_\cdot])$ with 
$\pcal=X$. To simplify the notation, we will drop the letter $N$ and consider $\xdot$ indifferently as a $\mathcal S
$-groupoid or a simplicial complex space, according to the needed properties. Notice that the stackification 
functor $\xdot\to\ccal\xdot$ corresponds to a fibrant resolution of simplicial presheaves. We conclude that
\begin{prop}\label{import1}
Let $p:X\to \ycal$ be an analytic stack. Then $p$ induces a groupoid equivalence $p:\xdot\to \ycal$ and a stack equivalence 
$\ccal\xdot\to \ycal$.
\end{prop}
An immediate consequence of this proposition is that, to work with a simplicial homotopy invariant, we can indifferently use any presentation and atlas of $\ycal$:  
\begin{cor}
Let $X,Z\to \ycal$ be two presentations of an analytic stack. Then $\xdot$ and $\zdot$ are equivalent 
groupoids and $\ccal\xdot$ and $\ccal\zdot$ are equivalent stacks. 
\end{cor}

In view of theorem \ref{fibranti}, given an atlas $X\to \ycal$, the stack $\ccal\xdot$ will simply be denoted
as $[X\times_\ycal X\rightrightarrows X]$ or $[X_1\rightrightarrows X]$.

\begin{digr}\label{moltipl}
{\em For future reference, we write in detail the product structure that $\xdot$ inherits from the presentation
$p:X\to \ycal$ of an algebraic stack. By definition of fibered product in the category of groupoids, an element
of the complex space $X_1=X\times_\ycal X$ is a triple $\uli{a}=(a_1, a_2, \alpha)$, where $a_1, a_2\in X$ 
and $\alpha:p(a_1)\stackrel{\cong}{\to} p(a_2)$ is an isomorphism. Notice that there is a canonical 
holomorphic map $e:X\to X_1$, the one sending $a\in X$ to $(a,a, {\rm id}_a)\in X_1$. Let $\partial_0$ and $\partial_1$   
 
\be
X\tms_\mathcal YX \stackrel{\d_0}{\underset{\d_1}{\rightrightarrows}} X
\ee
be the holomorphic  maps involved in the definition of fiber products, i.e. making the following diagram commutative
\be\label{str}
{\begin{aligned}
\xymatrix{X\tms_\mathcal YX\ar[r]^{\hskip 10pt \d_0}\ar[d]^{\d_1}& X\ar[d]^p\\
X\ar[r]^p&\ycal.}
\end{aligned}}
\ee
The multiplication $m$ is the holomorphic map associated to the pair $(\partial_0\circ pr_1, \partial_1\circ pr_2)$ and whose existence is consequence of the universality of fiber products:
\be\label{molt}
\xymatrix{X_1\times_X X_1\ar@{.>}[dr]|m\ar@/_/[ddr]_{\d_0\circ pr_1} \ar@/^/[rrd]^{\d_1\circ pr_2}  &&\\
& X\tms_\mathcal YX \ar[r]_{\;\;\d_1}\ar[d]^{\d_0} & X\ar[d]^p\\ 
& X\ar[r]^p&\ycal.}
\ee
$m$ is explicitely descibed as follows: if $a$ is as before and $b=(b_1,b_2,\beta)$ is another point of
$X_1$, then $m(a,b)=(a_1, b_2, \beta\circ \alpha)$. By the commutativity of the previous diagram, 
we have $(\partial_0\circ m)(a,b)=a_1$ and $(\partial_1\circ m)(a, b)=b_2$. 

If $\ycal$ is a  Deligne-Mumford analytic stack then $\p_0$ and $\_1$ are \etale; $m$ also is \etale because in the commutative diagram 
$$
\xymatrix{X_1\times_X X_1\ar[r]^{\;\;pr_1}\ar[d]^{pr_2} & X_1\ar[d]^{\d_0}\\X_1\ar[r]^{\d_1} & X}
$$
we have that $pr_i$ are \etale, hence $\partial_0\circ pr_1$ is \etale and from the commutativity of the 
diagram (\ref{molt}), we conclude that $m$ is \etale. 
}
\end{digr}


\section{Simplicial parabolic holotopy presheaves $\pis_i$}\label{hol-presh}

The notion of Brody hyperbolicity which we will introduce in Section \ref{iperbolicita} will be rephrased in terms of {\em simplicial parabolic holotopy presheaves}. It turns out that this is a very convenient way to work in practice with such seemingly abstract definitions.
Those holotopy presheaves completely determine the {\em global} simplicial homotopy class of a groupoid, 
and since local simplicial weak equivalences between locally fibrant simplicial presheaves coincide with global 
simplicial weak equivalences, holotopy presheaves completely determine equivalences of stacks. In other 
words, a functor $F:\xcal\to \ycal$ is an equivalence between stacks if and only if $F_*$ induces isomorphisms 
between all holotopy  presheaves of $\xcal$ and $\ycal$.

We recall that the {\em parabolic $n$-th dimensional circle} is the simplicial set ${\sf S}^1_s:=\Delta^1/\partial\Delta^1$, where $\Delta^1$ is the standard $1$-dimensional simplex, seen as constant presheaf in the analytic site. As usual, in what follows, we let $\wedge$ be the monoidal structure in $\Delta^{op}\pref$ and ${\sf S}^n_s:={\sf S}^1_s\wedge\stackrel{n}{\cdots} \wedge {\sf S}_s^1$.

\begin{defin}\label{holot}
\begin{itemize}
\item[1)] For any complex space $U$ and simplicial presheaf $\xcal$, we set
$$
\pis_i(\xcal,x)(U)\stackrel{def}{=}\homs_{\hsolop}({\sf S}^i_s\wedge U_+,(\xcal,x))
$$
\item[2)] For any groupoid $\gcal$ and $U\in\scal_T$ we set $\pi_i(\gcal,g)(U)=\pis_i(N\gcal,g)(U)$.
\end{itemize}
\end{defin}
\noindent By definition, the presheaves $\pis_i$, called {holotopy presheaves}, induce isomorphisms if applied to local and global weak equivalences
or groupoid equivalences. We know already how to compute most of these presheaves for groupoids $\gcal$: because of Theorem \ref{null}, $\pis_i(\gcal,g)$ are constant to $0$ for all $i$ greater or equal to $2$. 
In general, it is extremely hard to compute $\pis_i$ of groupoids if $i=0,1$. 
The following result shows that, in the case the groupoid is a stack, these presheaves are related to some combinatorial data of the stack:

\begin{prop}\label{boh2}
Let $\zcal$ be a locally fibrant simplicial presheaf. Denote $\partial_I$ the composition $\partial_{i_1}\circ\cdots\circ \partial_{i_n}$ for a multiindex $I=(i_1,\cdots,i_n)$, where $\partial_{i_j}:\zcal_j\to\zcal_{j-1}$ are face morphisms. Then, $\pis_n(\zcal,z)(U)$ is the quotient set $A_n(U)/\sim$ where
\begin{itemize}
\item[1)] $A_n(U)=\{f:U\to \zcal_n\,\text{\rm such that}\; \partial_I\circ f=\partial_J\circ f\;,\; \forall I\neq J\,\text{\rm of length}\;n\}$
\item[2)] $\sim$ is the equivalence relation generated by $f\sim g$ if it exists an 
$H:U\to \zcal_1$ such that $\partial_0\circ H=\partial_I\circ f$ e $\partial_1\circ H=\partial_I\circ g$.
\end{itemize}
\end{prop}
\demo
One can get the given description of $A_n$ by explicitely writing down the simplicial set ${\sf S}^n_s$ and 
considering simplicial morphisms from ${\sf S}^n_s$. The fact that $\pis_n$ is a quotient set of $A_n$ follows by the fibrancy assumption and the equivalence relation translate the conditions coming from the usual simplicial cylinder object to express simplicial homotopies between moprhisms. 
\enddemo

\begin{remark}
{\em In particular, for an analytic stack $\zcal$ 
$$
\pis_0(\zcal,z)(U)=\{f:U\to \zcal_0\}/\!\sim
$$
where $f\sim g$ if there exists a morphism $H:U\to \zcal_1$ such that $\partial_0\circ H=f$ and $\partial_1\circ H=g$. Defining $\pis_0$ this way for general groupoids, would result in a non (simplicial, local) homotopy invariant presheaf.
}
\end{remark}


\subsection{Descent data and simplicial parabolic holotopy presheaves}\label{hol-pref}

In the previous subsection we have observed that the only holotopy presheaves of a groupoid that are 
relevant are in degree zero and one. Moreover, if the groupoid is a stack $\ycal$, the sections of these 
presheaves are expressible in terms of sections in the presheaves {\sf Ob}$(\ycal)$ and ${\sf Mor}(\ycal)$, or, more 
precisely, of the sections of $(N\ycal)_0$ and $(N\ycal)_1$ since by 
Definition \ref{holot} $\pis_i(\ycal,x)(U)=\pis_i(N\ycal, x)(U)$ . 

\begin{prop}\label{import2}
Let $X\to \ycal$ be an analytic stack. Then 
$$\pis_i(\ycal,y)\cong\pis_i(\xdot,x)\cong\pis_i(\ccal\xdot,x)$$
\end{prop}

\demo
Use Proposition \ref{import1}.
\enddemo 

In view of these isomorphisms, for any simplicially homotopy invariant considerations about $\ycal$, we will replace $\ycal$ with $\ccal\xdot$ for some appropriate choice of atlas $X$ of $\ycal$. 
Because of the relevance of the concept in the sequel, we explicitely recall

\begin{defin}\label{discesa}
Let $\gcal$ be a groupoid, $U$ a complex space and $\ucal=\{U_i\}_i$ a covering of $U$ for the strong topology. Then,
\begin{itemize} 
\item[1)] A {\em descent datum} relative to $\ucal$ in $\gcal$ is a pair $((A_i), (h_{ij}))$, also denoted $(A_i, h_{ij})$ with: $A_i$ objects of 
$\gcal(U_i)$ and $h_{ij}:{A_i}|_{U_{ij}}\to {A_j}|_{U_{ij}}$ isomorphisms, called {\em transition morphisms} satisying the {\em cocycle condition}: $h_{jk}\circ h_{ij}=h_{ik}$ on $U_{ijk}$. As always, $U_{ij}$ and $U_{ijk}$ stand for the double and triple intersections of the indixed complex spaces.  
The set of descent data will be denoted by ${\sf Dis}_\gcal(\ucal)$.
\item[2)] A {\em descent data morphism} between $(A_i,h_{ij})$ and $(B_i,g_{ij})$ in $\gcal$ and relative to a covering $\ucal$ is a collection of isomorphisms $\{\phi_i: A_i\to B_i\}$ respecting the relation $g_{ij}\circ\phi_i=\phi_j\circ f_{ij}$ on $U_{ij}$ for all $i,j$. 
\end{itemize}
\end{defin}

\begin{remark}
{\em In a covering $\ucal$ associated to a descent datum we will possibly allow $U_i=U_j$ for $i\neq j$.  
}
\end{remark}

Given a complex space $U$, we will denote by Cov$\,U$ the set of all the countable, locally finite open coverings 
such that $U_i\Sbs U$ for all $i$. Any covering of $U$ 
can be refined to one in Cov$\,U$. This set is filtering with 
respect to the relation $\ucal\preceq\ucal'$ if $\ucal'=\{U_i\}_i$ is finer than $\ucal$ and $U_i'\Sbs U_{\tau(i)}$, 
if $\tau:\nbb\to \nbb$ is the refining function.

Consider a refinement $\ucal\preceq\ucal'$ of two coverings of $U$. There is a correspondence $r_{\ucal,\ucal'}:{\sf Dis}_\gcal(\ucal)\to {\sf Dis}_
\gcal(\ucal')$ which associates to a descent datum $(A_i,h_{ij})$ the datum $({A_i}|_{U_i'}, h'_{rs})$, where the new transition morphisms $h'_{rs}$ on the double intersections $U'_{rs}$ of complex spaces contained in $U_i$ are defined as the identity: ${A_i}|_{U'_{rs}}\to {A_i}|_{U'_{rs}}$. The same can be said for a descent data morphism $\phi=\{\phi_i\}_i$ between $a$ and $b$: it induces a morphism $r_{\ucal,\ucal'}(a)\to r_{\ucal,\ucal'}(b)$ in a unique way. Thus we have the following
\begin{lemma}\label{expl}
Let $X\to\ycal$ be an analytic stack. The stack $\ccal\xdot$ associated to the groupoid $\xdot$ is explicitely described by the following:
\begin{itemize}
\item[1)] ${\sf Ob}(\ccal\xdot)(U)=\varinjlim\limits_{U\in Cov\,U} {\sf Dis}_{\xdot}(\ucal)$;
\item[2)] let $\sf r$ and $\sf s$ be two objects rapresented by descent data $r$ and $s$, and defined on the same covering $\ucal$ (this is not restrictive). Then 
$$
\homs_{\ccal\xdot}({\rm r},{\rm s})=\big\{{\rm descent\> data\>(iso)morphisms}\>\> 
r\to s\big\}
$$
(see Definition \ref{discesa}).   
\end{itemize}
\end{lemma}  
\begin{remark}
{\em In \cite[Lemma 3.2]{lm} a different description of the stack associated to a prestack is given. In particular the objects are simply descent data unidentified in the direct limit. We think this is not a useful definition as the following example shows: let $\ycal=Y$ be a complex space and $X\to \ycal$ be a presentation with $X=\amalg_i B_i$, $B_i$ non trivial open subspaces. Then the canonical morphism $\ccal\xdot\to\ycal$ is a (presheaf) isomorphism, since the only isomorphisms in {\sf Mor} $\ccal\xdot$ are the identity morphisms. The objects of $\ccal\xdot$ over $U$ are the sections $\ycal(U)$, thus they must be all the holomorphic maps $U\to Y$. We do have a canonical surjective correspondence from the descent data over $U$ to holomorphic maps $U\to Y$, but this is not injective, factoring precisely through the relation defining the direct limit over the coverings of $U$. 
}
\end{remark}
The notion of descent data given in Definition \ref{discesa} may be expressed in terms of holomorphic maps if the groupoid 
$\gcal$ in question is $\ccal\xdot$. Given a covering $\ucal\in{\rm Cov}\,U$, $\ucal=\{U_i\}_i$ , a descent datum $r$ 
(on $U$) in $X$ relative to $\ucal$ is a pair $((r_i:U_i\to X)_i, (f_{ij}:U_{ij}\to X_1)_{ij})$ with $r_i$ and $f_{ij}$ holomorphic maps such that:
\bit
\item[$(\star)$]$\>\> r_{i|U_{ij}}={\d}_0\circ{ f}_{ij}, \>\>{r}_{j|U_{ij}}=\d_1\circ f_{ij}$ su $U_{ij}$;
\item[$({\star\star})$] $ f_{ij}$: $ m( f_{ij}\tms  f_{jk})= f_{ik}$ su $U_{ijk}$ (cocycle relation) 
\eit
and, like before,
$$
 \varinjlim\limits_{\ucal\in{\rm Cov}\,U}
\Bigl\{\bigl((r_i:U_i\to X)_i, (f_{ij}:U_{ij}\to X_1)_{ij}\bigr)\Big\}=\varinjlim\limits_{\ucal\in{\rm Cov}\,U} {\sf Dis}_{\xdot}(\ucal)={\sf Ob}(\ccal\xdot)(U)
$$

We are ready now to describe the zeroth and first holotopy presheaves of an analytic stack by means of the complex structure of any of its atlases:

\begin{theorem}\label{pi0}
Let  $p:X\to\ycal$  be an analytic stack and $U$ a complex space. Then
$$
\pi_0(\ycal,y)(U)\cong \varinjlim\limits_{\ucal\in{\rm Cov}\,U}
\Bigl\{\bigl((r_i:U_i\to X)_i, (f_{ij}:U_{ij}\to X_1)_{ij}\bigr)\Bigr\}/{\sim_0}
$$
where $\sim_0$ is the equivalence relation generated by $(r_i,f_{ij})\sim_0 (s_i,g_{ij})$ if and only if there exist holomorphic maps $\phi_i:U_i\to X_1$ such that
\begin{itemize}
\item[1)]\label{c1} $\partial_0\circ \phi_i=s_i$ e $\partial_1\circ \phi_i=r_i$ for all $i$;
\item[2)]\label{c2} $m(f_{ij}\times\phi_j)=m(\phi_i\times g_{ij})$. 
\end{itemize}
\end{theorem}\noindent
An isomorphism $\phi$ between two descent data $r=(r_i, f_{ij})$, $s=(s_i, g_{ij})$ is a collection of holomorphic  maps $\phi_i:U_i\to X_1$ such that 
\begin{itemize}
\item[$1'$)] $\partial_0\circ \phi_i=r_i$ and $\partial_1\circ \phi_i=s_i$ for all $i$;
\item[$2'$)] $m(f_{ij},\phi_j)=m(\phi_i, g_{ij})$ for all $i,j$.
\end{itemize}
Each collection $\{\phi_i\}_i$ determines a class in the filtered colimit, over the coverings $\ucal$ of $U$,
of isomorphisms between the descent data $r$ and $s$. Representatives of sections of $\pis_1(\ycal,y)(U)$ are
(classes of) automorphisms $(\phi)_i$ of $r$, for $r$ ranging in $\varinjlim\limits_{\ucal\in{\rm Cov}\,U}
{\sf Dis}_{\xdot}(\ucal)$. 
\begin{theorem}\label{pi1}
Let $p:X\to\ycal$ an analytic stack and $U$ be a complex space. Then,
$$
\pi_1(\ycal,y)(U)\cong \varinjlim\limits_{\ucal\in{\rm Cov}\,U}(\phi_i)_i/\!\!\sim_1
$$
where $\phi=(\phi_i)_i$ is an automorphism  of a descent datum $r$ in $\xdot$ relative to $\ucal$. 
If $\phi$ e $\psi$ are automorphisms of descent data $r$ and $s$, respectively, relative to $\ucal=\{U_i\}_i$, then $
\sim_1$ is the equivalence relation generated by $\phi\sim_1\psi$ if and only if there exists and isomorphism  
between $\partial_0\circ\phi=s$ and $\partial_0\circ\psi=r$, i.e. holomophic maps $H_i: U_i\to X\times_\ycal X$ such that $\partial_0\circ H_i=r_i$ and $\partial_1\circ H_i=s_i$ satisfying the second of the conditions listed in the Theorem \ref{pi0}.
\end{theorem}
\demo (of Theorems \ref{pi0} and \ref{pi1})
We use Proposition \ref{boh2}. By the Yoneda lemma the set $A_0$ coincides with the sections

$$
\ccal\xdot(U)={\sf Ob}(\ccal\xdot(U))=\varinjlim\limits_{\ucal\in{\rm Cov}\,U}{\sf Dis}_{\xdot}(\ucal).
$$
Thus, $\pis_0(\ccal\xdot,x)(U)$ is the quotient of this set by the equivalence relation given by the simplicial 
homotopy, that in turn is described in the part (2) of the above proposition. The same argument gives the 
statement  for the groups $\pis_1(\ccal\xdot,x)(U)$. Finally, we identify the presheaves $\pis_0(\ccal\xdot,x)$ 
and $\pis_1(\ccal\xdot,x)$ with those of $\ycal$ by means of the Proposition \ref{import2}.
\enddemo

\begin{defin}\label{sez-cost}
Let $\pcal$ be a presheaf on $\scal_T$. A section $\sigma\in \pcal(U)$ is {\em constant} if it lies in the image of the 
map $c^*:\pcal({\sf pt})\to \pcal(U)$, where $c:U\to {\sf pt}$.  
\end{defin}



\section{Hyperbolicity}\label{iperbolicita}

{\bf General assumption}: {\em For the time being, given an analytic stack $X\to \ycal$, we will assume that the two complex spaces $X$ and $X\times_\ycal X$ are {\em reduced}.}
\smallbreak

The classical Brody's Theorem claims that two notions of hyperbolicity for complex spaces coincide. One is rooted in metric aspects of 
the complex space, the other is defined in terms  of certain holomorphic maps. We are not aware 
of any possible candidates for the analogues of these notions about stacks in the literature. This section
is devoted to providing ours.


\subsection{Brody hyperbolicity}

In the paper \cite{BT2} we have given the following definition 

\begin{defin}\label{simpl-iper}
A simplicial presheaf $\pcal$ is {\em Brody hyperbolic} if
\begin{itemize}
\item[1)] is simplicially locally fibrant and 
\item[2)] the projection $p_{\xcal}:\C\times\xcal\to\xcal$ induces set bijections
\begin{equation}
\homs_{\hsolo}(\xcal, \pcal)\stackrel{\cong}{\to} \homs_{\hsolo}(\C\times\xcal,\pcal)
\end{equation} 
for all $\xcal\in\pref$. 
\end{itemize}
\end{defin}
Since a groupoid can be seen as a simplicial presheaf by means of the functor $N$, we will use the 
same definition

\begin{defin}\label{brody-iper1}
A groupoid $\gcal$ is {\em Brody hyperbolic} if $N\gcal$ is a Brody hyperbolic simplicial presheaf. 
\end{defin}
Notice that a Brody hyperbolic groupoid is necessarily a stack. This definition can be rephrased in 
terms of holotopy presheaves:

\begin{prop}\label{rephr}
Let $\pcal$ a locally fibrant simplicial presheaf. The following conditions are equivalent:
\begin{itemize} 
\item[1)]\label{morvoe} $\pcal$ is Brody hyperbolic;
\item[2)]\label{nuova} $p_{\xcal}^* :{\sf Map}(\xcal,\pcal)\to {\sf Map}(\C\times\xcal,\pcal)$ is a weak equivalence of simplicial sets for any $\xcal\in\Delta^{op}\pref$, where ${\sf Map}$ is the simplicial mapping space;
\item[3)]\label{quella-imp} the simplicial holotopy presheaves are Brody hyperbolic, i.e. the projection $p_U:\C\times U\to U$ induces isomorphisms $p_U^*:\pis_i(\pcal,y)(U)\stackrel{\cong}{\to}\pis_i(\pcal,y)(\C\times U)$ for each $i$ and complex space $U$. 
\end{itemize}
\end{prop}
\demo 
2)$\Rightarrow$ 1). $\xcal$ is locally fibrant by assumption and 
 $$
 {\sf Hom}_{\hcal_s}(\xcal,\pcal)=\pi_0({\sf Map}(\xcal,\pcal))\cong \pi_0({\sf Map}(\C\times\xcal,\pcal))={\sf Hom}_{\hcal_s}
(\C\times\xcal,\pcal).
$$
 
1)$\Rightarrow$2). We show that the condition 2) is equivalent to $\pcal(U)
\cong \pcal(\C\times U)$ being a simplicial sets weak equivalence for all complex spaces $U$. By a small 
object argument there is a weak equivalence $\phi: \zcal_n\to \xcal$, where $\zcal_n$ are direct sums 
of simplicial sheaves of the kind $U\times \Delta^n$ for complex spaces $U$, because for any presheaf
there is a surjection onto it from a direct sum of complex spaces. Moreover, for any simplicial presheaf $
\tcal_\cdot$, the canonical morphism hocolim $\tcal_n\to \tcal_\cdot$ is a simplicial local weak equivalence.  
Since both domain and codomain are cofibrant objects (all objects are for the Joyal simplicial, 
injective, local model structure), by \cite[Corollary 9.7.5.\,\!(2)]{hirsc}, $\phi^*:{\sf Map}(\xcal,\pcal)\to 
{\sf Map}(\zcal,\pcal)$ is a simplicial weak equivalence. In turn this simplicial set is weakly equivalent to
$$
{\sf Map}(\ocolim_n (\amalg_U U\times \Delta^n), \pcal)\cong \olim_n {\sf Map}(\amalg_U U\times\Delta^n,\pcal)=\olim_{n,U}{\sf Map}(U,\pcal).
$$ 
Similarly, 
$$
{\sf Map}(\C\times\xcal,\pcal)\cong \olim_{n,U}{\sf Map}(\C\times U,\pcal)
$$
since finite limits commute with filtered colimits. Thus, if we knew that
$$
\pcal(U)={\sf Map}(U,\pcal)\cong{\sf Map}(\C\times U,\pcal)=\pcal(\C\times U)
$$ 
for all complex spaces $U$, we would have that the homotopy limits are weakly equivalent and the condition 
2) is verified. To prove that the condition 1) implies that the simplicial sets
$\pcal(U)$ are weakly equivalent to $\pcal(\C\times U)$ for each $U$, we consider $\xcal={\sf S}^n_s\wedge 
U_+$ and use the pointed version of the condition 1). We conclude that 
$$
{\sf Hom}_{\hcal_s}({\sf S}^n_s\times U,\pcal)={\sf Hom}_{\hcal_s}({\sf S}^n_s\times
\C\times U,\pcal).
$$
 On the other hand, by adjunction, we have
\beq
{\sf Hom}_{{\hcal_s}_\bullet}({\sf S}^n_s\wedge  U,(\pcal,y))&=&{\sf Hom}_{{\hcal_s}_\bullet}({\sf S}^n_s, {\sf Map}_\bullet(U_+ ,(\pcal,y)))=\\ 
&&{\sf Hom}_{{\hcal_s}_\bullet}({\sf S}^n_s, (\pcal(U),y))\stackrel{def}{=}\pis_n(\pcal(U),y).\nonumber
\eeq 
Similarly, we prove that 
$$
{\sf Hom}_{{\hcal_s}_\bullet}({\sf S}^n_s\wedge U_+\wedge \cbb_+,(\pcal,y))=\pis_n(\pcal(\C\times U),y),
$$
hence the simplicial sets $\pcal(U)$ and $\ycal(\C\times U)$ are weakly equivalent.

We have yet to prove that one of the conditions 1) or 2) is equivalent to 
the condition 3). For $\pcal$ this is equivalent to have the sets (groups, when applicable)
 
\beq
\pis_i(\pcal,f)(U)&=&{\sf Hom}_{{\hcal_s}_\bullet}({\sf S}^i_s\wedge U_+,(\pcal,y))=\\
&&{\sf Hom}_{{\hcal_s}_\bullet}({\sf S}^i_s, {\sf Map}_\bullet(U_+,(\pcal,y)))=\pis_i(\pcal(U),y)\nonumber
\eeq
isomorphic to $\pis_i(\pcal(\C\times U),y)$, i.e. $p_U^*:\pcal(U)\to \pcal(\C\times U)$ is a weak homotopy equivalence, which is equivalent to condition 2) for what it has been previously proved.
\enddemo

The notion of Brody hyperbolicity we will more frequently use is the 3) of the previous 
proposition. 
\begin{defin}\label{br-campi}
\begin{itemize}
\item[1)] A presheaf  $\pcal$ is {\em Brody hyperbolic} if $p_U:\C\times U\to U$ induces bijections
$p_U^*:\pcal(U)\to \pcal(\C\times U)$ for any complex space $U$.
\item[2)]\!\! A groupoid $\gcal$ is\!\ {\em Brody hyperbolic} if the\! holotopy presheaves $\pis_i(\gcal,g)$ ( 
Defi-nition \ref{holot}) are hyperbolic for all $i$, hence only for $i=0,1$, because of Theorem \ref{null}.
\end{itemize}
\end{defin}


\subsection{Kobayashi hyperbolicity}

In the previous subsection we defined Brody hyperbolicity of a groupoid by first giving the same notion 
for a presheaf and then imposing that specification to the holotopy presheaves of the groupoid. The holotopy 
presheaves determine whether a groupoid is {\em Kobayashi hyperbolic}, as well. Classically, complex 
spaces Kobayashi hyperbolicity is a notion arising in the attempt to give complex spaces
a biholomorphically invariant distance. In general the best that can be done is endowing complex spaces of 
a biholomorphically {\em pseudodistance}. When on a complex space $X$ this happens to be a distance, $X$
is said to be {\em Kobayashi hyperbolic}.

The notion of Kobayashi hyperbolicity for groupoids is based upon the concept of relative analytic 
disc. 

\subsubsection{Discs and analytic chains}\label{dischi}
Let $\dbb$ be the unitary open disc in $\cbb$. We recall that for a complex space $U$, we have denoted $\text{Cov}\,U$ the set of countable, locally finite, open coverings $\ucal=\{U_i\}_i$ of $U$ such that $U_i\Sbs U$ for all $i$. If $\ucal\in \text{Cov}\, U$ and 
$\dcal=\{D_a\}_a\in \text{Cov}\, \dbb$ let $\{\dcal\times\ucal\}$ be the covering $\{D_a\times U_i\}_{ai}$ of 
$\dbb\times U$. The set $\text{Cov}\,\dbb\times \text{Cov}\, U$ is filtering in $\text{Cov}\,\dbb\times U$.

 Let $\pcal$ be a presheaf. A {\em relative analytic disc} of $\pcal$ on a compex space $U$ is 
an object of ${\sf F}\in\pcal(\dbb\times U)$. For any $z\in \dbb$, the same letter will refer to the inclusion $\{z\}\times U\hookrightarrow \dbb\times U$.
Let $r,s\in \pcal(U)$ be two sections and suppose there exists a relative analytic disc ${\sf F}$ and
two points $z_1, z_2\in \dbb$ such that $z_1^*{\sf F}=r$ and $z_2^*{\sf F}=s$. The sections $r$ and $s$ 
are then said to be {\em connected} by ${\sf F}$. 
A {\em relative analytic chain} on $U$ connecting $r$ to $s$ is the set ${\sf C}_{(r,s)}$ of the following data:
\begin{itemize}
\item[1)] a collection $r_0=r,\cdots ,r_k=s$ of sections;
\item[2)] $2k$ points $a_1,b_1,\cdots, a_k, b_k$ in $\dbb$;
\item[3)] $k$ relative analytic discs ${\sf F}_1,\cdots , {\sf F}_k$ such that the analytic disc ${\sf F}_i$ connects 
the sections $r_{i-1}$ and $r_i$, i.e. $a_i^*{\sf F}_i=r_{i-1}$ and $b_i^*{\sf F}_i=r_i$ for all $1\leq i\leq k$.
\end{itemize}

If a relative analytic chain ${\sf C}_{(r,s)}$ connects the section $r$ with the section $s$,
we call the pair $(r,s)$  {\em admissible}. If $\pcal=Y$ is a complex space, admissibility of all section
pairs in $\pcal({\sf pt})$ is equivalent to the topological connectedness of $Y$.

Endow the unitary open disc $\dbb$ of the Poincar\cce metric 
$$
ds^2=\dfrac{1}{(1-\vert z\vert^2)^2}dz\otimes d\bar z
$$
and denote with  $\r_{_{\D}}(p,q)$ the associated distance function between two points $a$ and $b$ in $\dbb$.
Then, for every chain ${\sf C}_{(r,s)}$ the nonnegative number 

\be\label{lcat}
l({\sf C}_{(r,s)})=\r_{_\D}(a_1,b_1)+\cdots+\r_{_\D}(a_k,b_k)
\ee
is, by definition, the {\em (Kobayashi) length} of the relative analytic chain ${\sf C}_{(r,s)}$. 
If $(r,s)$ is an admissible pair of sections, the nonnegative real number

\be\label{ddisc}
{\rm d}_{\rm Kob}^{\pcal}(r,s)=\inf\limits_{{\sf C}_{(r,s)}}l({\sf C}_{(r,s)})
\ee
defines a {\em pseudodistance} function on all the admissible pairs of sections in $\pcal(U)$ for all complex 
spaces $U$, called {\em Kobayashi pseudodistance} of $\pcal$. 

Arguing as in the case of complex spaces it immediately seen that morphisms of presheaves decrease the Kobayashi pseudodistance.

\begin{defin}
A presheaf  $\pcal$ is said to be {\em Kobayashi hyperbolic} if its Kobayashi pseudodistance is indeed a 
distance, hence if and only if ${\rm d}_{Kob}^\pcal(r,s)\neq 0$ for all admissible pairs $(r,s)\in \pcal(U)$ with $r\neq s$ and all complex spaces $U$.
\end{defin}

The Kobayashi hyperbolicity for a groupoid is defined as follows:

\begin{defin}\label{kob-iper}
A groupoid $\gcal$ is {\em Kobayashi hyperbolic} if the presheaves $\pis_0(\gcal,g)$ and $
\pis_1(\gcal,g)$ are Kobayashi hyperbolic. 
\end{defin}

For general groupoids, we cannot go much further in expliciting what this amounts to; however, if $\gcal$ happens to be an analytic stack, the description of holotopy presheaves given in the Section \ref{hol-presh} may be used to express the Kobayashi hyperbolicity in terms of holomorphic maps between complex spaces. 

Thus, consider an analytic stack $p:X\to \ycal$. Since Kobayashi hyperbolicity of $\ycal$ is a prescription
on the holotopy sheaves, we can replace $\ycal$ with the stack $\ccal\xdot$, because of Proposition \ref{import2} and, by Theorem \ref{pi0}, this condition can be rephrased using descent data in $X$. 
In this case, a relative analytic disc is determined by a descent datum $\sf F$ on $\dbb\times U$ and therefore by an open covering $\dcal\times \ucal$ of $\D\tms U$ and by holomorphic maps

$$
F_{ai}:D_a\times U_i\to X, \; F_{ai,bj}:D_{ab}\times U_{ij}\to X_1=X\times_\ycal X,
$$
where we wrote $D_{ab}\times U_{ij}$ for the set $(D_a\times U_i)\cap (D_b\times U_j)$. These maps
satisfy the conditions $(\star)$ and $(\star\star)$ mentioned in the Subsection \ref{hol-pref}. Keeping the usual 
notation, this will be denoted as
\be\label{disc1}
{\sf F}={\sf F}_{\dcal\tms\ucal}=\big(F_{ai}:D_a\tms U_i\to \D\tms U , F_{aibj}:D_{ab}\tms U_{ij}\to X_1\big)
\ee
with the coverings morphisms ranging among those of the kind $\amalg_{ai}(D_a\tms U_i)\to\D\tms U$. 

The covering $\{D_a\times U_i\}_{ai}$ induces on $\{z\}\times U$ an open covering $\{U_{ai}\}_{ai}$ of $\{z\}
\times U$ by setting $U_{ai}=(D_a\tms U_i)\cap  \{z\}\times U$ and the holomorphic maps $F_{ai}$ and 
$F_{aibj}$ restrict to holomorphic maps on $U_{ai}$ and their intersections, satisfying  $(\star)$ and $
(\star\star)$. Notice that among the open subspaces $\{U_{ai}\}$ there may be some indexed by different 
pairs but are coincident are coincident as subspaces. A descent data $r=(r_i, f_{ij})$ on $U$ relative to a covering $\{U_i\}_i$ determines a descent data $r$ relative to the covering $\{U_{ai}\}_{ai}$ by letting $r_{ai}=r_i$ e $f_{aibj}=f_{ij}$. 

\begin{remark}
{\em
\begin{enumerate}
\item[1)] Let $f:\xcal\to\ycal$ be a $1$-morphism of groupoids. Then for any complex space $U$ and pair 
of admissible sections $r,s\in \pis_i(\xcal,x)(U)$, the sections $f_*(r)$ and $f_*(s)$ are admissible in 
$\pis_i(\ycal,y)(U)$ and 
$$
{\rm d}_{\rm Kob}^\ycal(f_*(r),f_*(s))\leq {\rm d}_{\rm Kob}^\ycal(r,s).
$$
It follows that ${\rm d}_{\rm Kob}$ is a simplicial, local weak equivalence invariant of stacks.\medskip
\item[2)] If $\ycal$ is a complex space $Y$, then $\pis_i=0$ for $i>0$ and $\pis_0(\ycal,y)\cong Y$ as (pre)sheaf.
The definitions of Brody and Kobayashi hyperbolicity we have given involve all sections, holomorphic 
 maps in this case, $U\to Y$ as opposed to just points $y\in Y$ as it is classically the case. It is easy to see 
that the relevant notions of hyperbolicity coincide.\medskip 
\item[3)] At the other extreme end, let $\ycal=\bcal G$ be the classifying stack of a Lie group $G$. Then $X={\sf pt}$
and $G$ has to be finite, if we wish to restrict to Deligne-Mumford analytic stacks. Then  there are no admissible
pairs of sections in $\pis_0$ and $\pis_1$, being similar to the case of a discrete complex space.   
\end{enumerate}
}
\end{remark}


\section{The coarse moduli space}\label{aus} 

Although the set $\ccal\xdot(U)=\varinjlim\limits_{\ucal\in{\rm Cov}\,U} {\sf Dis}_{\xdot}(\ucal)$ (see 
Lemma \ref{expl}) is closely tied to the complex structure of an atlas $X$ of an analytic stack $\ycal$, it is unclear how 
to metrize it in a usable way and arguments employing complex variables theory, such as those necessary to prove 
Brody theorem, seem not possible. Alternatively, it is possbile to link in a natural way this set (see to the proof of Lemma \ref{fun}) to a complex space, the coarse moduli space of $\ycal$, which we will denote by $Q(\ycal)$.


\subsection{Compactness}

In the classical Brody Theorem compactness is essential; similarly here we will need some finiteness condition
on the stacks in order to proceed further. 

Given an analytic stack with presentation $p_X:X\to\ycal$, we will denote by $j_X$ the holomorphic map
$$
j_X=(\partial_0, \partial_1):X_1=X\times_\ycal X\to X\times X.
$$
Let $\ycal$ be an analytic stack. An {\em open substack} $\ocal\subset\ycal$ is a full category of $\ycal$ such that
\begin{itemize}
\item[1)] for any object $x\in\ocal$ all objects in $\ycal$ isomorphic to $x$ are also in $\ocal$;
\item[2)] the inclusion morphism $\ocal\subset\ycal$ is represented by open immersions.
\end{itemize}
If every covering $\amalg_{\a}\ocal^\a\to \ycal$ by open substacks admits a finite subcovering and for one presentation $W\to \ycal$ the map $j_W$ is proper, we will say $\ycal$ is {\em compact by open coverings}.

The following three properties will serve to our purpose:

\begin{prop}\label{compat}
Let $X\to \ycal$ be an analytic stack.  The following are equivalent:
\begin{enumerate}
\item[1)] $\ycal$ is compact by open coverings;
\item[2)] there exist two presentations $p_X:X=\amalg_{i=1}^N X^{(i)}\to 
\ycal$, $p_Z:Z=\amalg_{i=1}^M Z^{(k)}\to \ycal$ with $X^{(i)} $, $i=1,\ldots,N$, $Z^{(k)}$, $k=1,\ldots,M$ connected  and  a holomorphic map $\phi:X\to Z$ over $\ycal$ such that $\phi|_{X^{(i)}}$ is an open embedding $X^{(i)}\hookrightarrow  Z^{(k)}$ for some $k$ with $\phi(X^{(i)})\Sbs Z^{(k)}$ and the map $j_Z:Z\times_\ycal Z\to Z
\times Z$ is proper.
\end{enumerate}
Under the equivalent conditions $\rm 1), 2)$ we will say that $p_Z:Z\to \ycal$ and $p_X:X\to \ycal$ are {\em adapted} presentations.
\end{prop}
\demo
Let $p_X:X=\amalg_{i=1}^N X^{(i)}\to \ycal$, $p_Z:Z=
\amalg_{i=1}^N Z^{(i)}\to \ycal$ be adapted presentations $\ycal$ and $\amalg_{\a}\ocal_\a\to \ycal$ an open covering by substacks. 
Then $Z\tms_\ycal\ocal_\a\hookrightarrow Z$ is an open immersion of complex spaces for every $\a$. It follows that $\amalg_{\a}(Z
\tms_\ycal\ocal_\a)$ restricted to $X$ is an open covering of $X$.  If $\amalg_{j=1}^n(Z\tms_\ycal\ocal_{\a_j})$ covers 
$\overline X$ then $\amalg_{j=1}^n\ocal_{\a_j}\to \ycal$ is a covering of $\ycal$.

Conversely, suppose that $\ycal$ is compact by open coverings and let $p_W:\amalg_{i=1}^{+\infty}W^{(i)}\to\ycal$ a 
presentation of $\ycal$ with $p_W$ open and $j_W:W\times_\ycal W\to W\times W$ proper. Cover each $W^{(i)}$ by relatively compact balls and apply the shrinking lemma to obtain two presentations $p_{\hat{X}}:\hat{X}=\amalg_{i=1}^{+\infty} \hat{X}^{(i)}\to \ycal$, 
$p_{\hat{Z}}:\hat{Z}=\amalg_{i=1}^{+\infty} \hat{Z}^{(i)}\to \ycal$ with $\hat{X}^{(i)}\Subset \hat{Z}^{(i)}$, $i=1,\ldots,N$. Consider the $\hat{X}_1=\hat{X}\times_\ycal \hat{X}$  saturations of the open spaces $\hat{X}^{(i)}$ (i.e. $\partial_1\partial_0^{-1}(\hat{X}^{(i)})$) and call them $X^{(i)}$. Let $\ocal_i\subset \ycal$ be the open substacks
corresponding to them. Only finitely many of them, say $\{\ocal_{i_s}: j=1,\cdots, N\}$ will be necessary to 
cover $\ycal$, by assumption. Thus $X:=\amalg_{s=1}^N X^{(i_s)}\to \ycal$ is an atlas. Doing the same with $\hat{Z}$
we produce an atlas $Z:=\amalg_{k=1}^M Z^{i_k}$ of $\ycal$. Since $j_W$ is proper, $X^{(i_s)}\Sbs Z^{(i_k)}$ for 
some $k$ showing that $X$ and $Z$ are the sought for atlases.
 
This proves that conditions {\rm 1), 2)} are equivalent. 
 \enddemo



\begin{defin}
A stack satisfying either of the conditions of the proposition \ref{compat} will be called {\em compact}.
\end{defin}

\begin{remark}\label{procpt}
 {\em Facts:
\begin{itemize}
\item[1)] if a stack admits a compact atlas, the stack is compact;
\item[2)] compactness is invariant by equivalence of stacks;
\item[3)] if $\ycal$ is a compact stack the diagonal morphism $\Delta:\ycal\to\ycal\tms\ycal$ is proper.
\end{itemize}
}
\end{remark}
1) and 2) are clear. In order to prove 3) fix two adapted presentations $p_Z:Z\to \ycal$ and $p_X:X\to \ycal$ and consider a diagram
$$
{\begin{split}
\xymatrix{S\tms_{(\ycal\tms\ycal)}\ycal\ar[r]^{}\ar[d]^{}& \ycal\ar[d]^\Delta\\
S\ar[r]^{} & \ycal\tms\ycal}
\end{split}}
$$
where $S$ is a complex space. The left vertical line is a morphism of complex spaces $\Delta$ being representable. 
Let $\{\xi_n\}$ be a sequence in the complex space $S\tms_{(\ycal\tms\ycal)}\ycal$: each $\xi_n$ corresponds to an 
element $[s_n,y_n,(a(s_n),b(s_n)),(\phi'_n,\phi''_n)]$ of the stack $S\tms_{(\ycal\tms\ycal)}\ycal$ where $s_n\in S$, 
$y_n\in{\sf Ob}(\ycal)$, $(a(s_n),b(s_n))\in{\sf Ob}(\ycal\tms\ycal)$, and $(\phi'_m,\phi''_m)$ is an isomorphism $
(a(s_n),b(s_n)){\cong}(y_n,y_n)$. Suppose, by contradiction, that $\{\xi_n\}$ has no accumulation points whereas a 
subsequence $\{s_m\}\sbs\{s_n\}$ converges to a point $\underline s\in S$. Let $\{x_m\}$, $\{x'_m\}$, $\{x''_m\}$ be 
sequences in $X$ satisfying $p_X(x_m)\cong y_m$, $\phi'_m:p_X(x'_m)\stackrel{\cong}{\to} a(s_m)$, 
$\phi''_m:p_X(x'_m)\stackrel{\cong}{\to} b(s_m)$. 
In view of 
compactness, we may assume that there exists $z,z',z''\in Z$ and $\phi',\phi''\in Z_1$ such that $x_m\to z$, $x'_m\to 
z'$, $\phi'_m\to\phi'$, $\phi''_m\to\phi''$. The sequence $[s_m,x_m,(a(s_m),b(s_m)),(\phi'_m, \phi''_m)]$ 
in the complex space $S\times_{\ycal\times\ycal} Z$ converges to $[\underline s,z,(a(\underline s), b(\underline s)),
(\phi',\phi'')]$ in $S\times Z\times Z_1\times Z_1$. The latter point lies, infact, in $S\times_{\ycal\times\ycal} Z$, because $(\phi', \phi'')$ is an isomorphism between $(a(\underline s), b(\underline s))$ and $\Delta(p_Z(z))$ (use
the morphism $j$). By applying 
$$
\mathrm{id}_S\times_{\ycal\times\ycal} p_Z: S\times_{\ycal\times\ycal} Z\to 
S\times_{\ycal\times\ycal} \ycal
$$
we deduce  that 
$$[s_m,y_m,(a(s_m),b(s_m)),(\phi'_m, \phi''_m)]
\to [\underline s,p_Z(z),(p_Z(z'),p_Z(z'')),(\phi',\phi'')].
$$

 
\subsection{Existence of the complex structure}\label{exist}

In the paper \cite{grezzo} we have proved that to a flat analytic groupoid 
$\xcal=\{s, t: R\rightrightarrows X\}$ with  $j=(s, t):R\to X\times X$ finite ( or equivalently $j$ closed and $\xcal$ with finite stabilizer) is associated a GC quotient. This is a complex space $Q(\xcal)$ with an holomorphic map 
$q: X\to X/R\to Q(\xcal)$ satisfying a number of conditions. Here $X/R$ is the quotient sheaf associated with the 
equivalence relation induced by $R$. That result applies to the analytic groupoid $\xcal=\big\{X\tms_\mathcal YX
\stackrel{\d_0}{\underset{\d_1}{\rightrightarrows}}X\big\}$ where $p_X:X=\amalg_{i=1}^N X^{(i)}\to \ycal$ is a 
presentation of the compact  Deligne-Mumford analytic stack $\ycal$. 


We can state another property of a stack equivalent to compactness:

\begin{prop}\label{cptdiag}
Let $p:X\to\ycal$ be an analytic stack with $p$ open and admitting a coarse moduli space $Q(\ycal)$. Then $\ycal$ is compact if and only if $Q(\ycal)$ is a compact complex space and $j$ is proper for some presentation. In that case the diagonal morphism $\Delta_Q:\ycal\to\ycal\tms_{Q(\ycal)}\ycal$ is proper.
\end{prop}

{\em Proof.} If $\ycal$ is compact, $Q(\ycal)$ is compact as well. Conversely,  let $W\to \ycal$ be a presentation.
As in the proof of proposition \ref{compat} we can always produce two atlases $A=\amalg_{i=1}^\infty A^{(i)}$ 
and $B=\amalg_{j=1}^\infty B^{(j)}$ satisfying 
all the conditions for the compactness of $\ycal$ except the finiteness of the number of the connected components.

Let $q_A:A\to Q(\ycal)$ and $q_B:B\to Q(\ycal)$ be the canonical holomorphic maps and $\{U_k: k=1,\cdots,N\}$
a finite subcover of $\{q_A(A^{(i)}), q_B(B^{(j)})\}_{i,j}$. Let $X^{(k)}$ be the $A^{-1}$-saturated of the open set 
$q_A^{-1}U_k$ and $Z^{(k)}$ be the $B_1$-saturated of $q_B^{-1} U_k$.
Then $\amalg_{k=1}^N X^{(k)}$ and $\amalg_{k=1}^NZ^{(k)}$ are two atlases with the required properties. 
Indeed, each $X^{(k)}$ and $Z^{(k)}$ are the atlas of an open substack of $\ycal$ and the union of them is 
surjective on $\ycal$.

In order to prove that $\Delta_Q$ is proper set $Q=Q(\ycal)$ and consider the following diagram

$$
\xymatrix{S\tms_{(\ycal\tms_{Q}\ycal)}\ycal \ar[r]^{}\ar[d]_{} &\ycal\ar[r]^{\rm id}\ar[d]_{\Delta_Q} & \ycal\ar[d]_{\Delta} \\ S\ar[r]^{} & \ycal\tms_{Q}\ycal \ar[r]^{} & \ycal\tms\ycal}
$$
where $S$ is a complex space. Then one checks that right side square is cartesian i.e. $\ycal$ is equivalent to the fiber product of $\ycal$ and $\ycal\tms_{Q}\ycal $ over $ \ycal\tms\ycal$; thus properness of $\Delta_Q$ follows from that of $\Delta$ (see Remark \ref{procpt}).
\enddemo

\section{Topological and metric structures}
\subsection{Distances}\label{top}
For the time being, we will only consider compact,  Deligne-Mumford analytic stacks.
Let $\ycal$ one such analytic stack with atlases $X\Sbs Z$, $X=\amalg_{i=1}^NX^{(i)}$, $Z=\amalg_{i=1}^NZ^{(i)}$ (see Section \ref{aus}). We can assume $X_i$ and $Z_i$ are Stein.

Fix a distance ${\rm d}:Q(\ycal)\tms Q(\ycal)\to\R_{\ge 0}$ on the quotient space $Q(\ycal)$ determined by a differentiable length function $H$. As we have previously 
introduced, $q=q_X$ is the projection $X\to Q(\ycal)$. Even if $q^*H$ is only a pseudolength function neverthless to $q^*H$ is associated a  distance on any connected component of $X$: this because $q$ is locally proper and equifinite fibres. This is the same distance induced by the restriction of $q_Z^*H$ to $X$. These distances on the connected components can be assembled together to a unique distance ${\rm d}_X$
on all $X$ in the following way.
Fix two points  $x\neq y\in X$; a (piecewise differentiable) path $\gamma$ through $x$ and $y$ is a set
$\{\gamma_1,\gamma_2,\cdots, \gamma_m\}$ of paths $[0,1]\to X$ such that:
\begin{enumerate}
\item[1)] $\gamma_1(0)=x$ and $\gamma_m(1)=y$;
\item[2)]  $q(\gamma_{i+1}(0))=q(\gamma_i(1))$ for all $1\leq i\leq m+1$.
\end{enumerate}
If $l(\gamma_i)$ denotes the length of $\gamma_i$ with respect to the distance on the connected component of $X$ in which the image of $\gamma_i$ lies, the positive real number $l(\gamma)=\sum_{i=1}^ml(\gamma_i)$ is the length of $\gamma$ by definition. We then set
$$
{\rm d}_X(x,y)=\inf\limits_{\g}l(\g).
$$
We proceed likewise  for the complex space $X_1=X\times_\ycal X$, considering on $X_1$ the length 
functions $(q\circ\partial_0)^*H$, $(q\circ \partial_1)^*H$ which give rise to the same distance making 
$\partial_0$ and $\partial_1$ into local isometries, since $q\circ \partial_0=q\circ \partial_1$. A similar argument applies to $X_2:=X_1\times_{\partial_0,X,\partial_1}X_1$ and to the multplication $m:X_2\to X_1$
that becomes a local isometry, as well.

The distance just introduced allow to metrize in a natural way the sets ${\sf Dis}_{\xdot}(U)$, $\pi_0(\ycal,y)(U)$ e $\pi_1(\ycal,y)(U)$ as follows. Let $r=({r}_i,{f}_{ij})\neq s=(s_i,g_{ij})$ two descent  data in ${\sf Dis}_{\xdot}(\ucal)$, with $\ucal\in \text{Cov}\,U$; define
 
\be\label{dist}
\delta_U(r,s)=\sup\limits_{\stackrel{u\in U_i}{i\in\N}}\bigl\{{\rm d}_X(r_i(u),s_i(u))\bigr\}+\sup\limits_{\stackrel{u\in U_{ij}}{i,j\in\N}}\bigl\{{\rm d}_{X_1}(f_{ij}(u),g_{ij}(u))\bigr\}.
\ee
The distance function $\delta_U(r,s)$ is invariant by restriction, that is $\delta_U(r,s)=\delta_U(r_{|\ucal'},s_{|\ucal'})$ 
for any $\ucal'\succeq\ucal$, thus it is defined for pairs of objects in ${\sf Ob}(\ccal\xdot (U))$. In turn it induces a distance function $\delta{\pi_0,U}$ on $\pis_0(\ycal,y)(U)$: for $\alpha,\beta\in \pis_0(\ycal,y)(U)$ 
\be\label{distp0}
\delta_{\pi_0,U}(\a,\b)=\inf\limits_{r\in\a,s\in\b}\delta_U(r,s).
\ee
For $\pis_1(\ycal,y)(U)$ we proceed similarly. Given $\ucal=\{U_i\}_i\in \text{Cov}\,U$, the elements of $\pis_1(\ycal,y)(U)$ are represented by pairs $[r,\phi]$ of a descent datum $r$ and $\phi=\{\phi_i\}_i$ is an
automorphism of $r$ (see Theorem \ref{pi1}). The distance $\delta$ between two such pairs $(r, \phi)$ and 
$(s,\psi)$ is defined as
  
\be\label{1dist}
\delta'_U\bigl((r,\phi), (s,\psi)\bigr)=\delta(r,s)+\sup\limits_{\stackrel{u\in U_i}{i\in\N}}\bigl\{{\rm d}_{X_1}\bigl(\phi_i(u),\psi_i(u)\bigr)\bigr\}
\ee
and the distance between two classes $\alpha,\beta\in \pis_1(\ycal,y)(U)$ is

\be\label{distp1}
\delta_{\pi_1,U}(\a,\b)=\inf\limits_{\stackrel{(r,\phi)\in\a}{(s,\psi)\in\b}}\delta'_U\bigl((r,\phi),(s,\psi)\bigr).
\ee
\begin{prop}\label{distp2}
$\delta_{\pi_0,U}$ e $\delta_{\pi_1,U}$ are distances.
\end{prop}
\demo
Let us introduce the following notation: $\widetilde{r}\neq\widetilde{s}$ are two classes in $\pis_0(\ycal,y)(U)$ for a fixed 
complex space $U$ and descent data $r=(r_i,f_{ij})$ and $s=(s_i,g_{ij})$ in $X$. It suffices to show that $
\delta_{\pi_0,U}(\widetilde{r}, \widetilde{s})>0$, since the other axioms of being a distance are readily seen to be 
satisfied. Recall that, given a descent datum $r$, the composition $q\circ r$ denotes the well defined 
holomorphic map $U\to Q(\ycal)$ defined locally by $q\circ r_i$ and gluing all these  maps over $U$. 
We split the proof in two complementary cases: {\em $q\circ r\neq q\circ s$}, {\em $q\circ r=q\circ s$.} 
\medbreak
\noindent{\em $q\circ r\neq q\circ s$.} Under this assumption we have that $\sup\limits_{u\in U_i}\left\{{\rm d}_X
\left(r_i(u),s_i(u)\right)\right\}\geq N>0$: let $V\subset U$ be an open set such that $\sup\limits_{v\in V}
\left\{{\rm d}_{Q(\ycal)}(q\circ r(v), q\circ s(v))\right\}=N>0$. Then 
$$
\delta_U(r,s)\geq \sup\limits_{\stackrel{u\in U_i}{i\in\N}}\left\{{\rm d}_X\left(r_i(u),s_i(u)\right)\right\}\geq N>0
$$ 
since the length function on $X$ inducing $d_X$ is $q^*H$, where $H$ is a length function on $Q(\ycal)$,
hence d$_X(x_1,x_2)\geq {\rm d}_{Q(\ycal)}(q(x_1),q(x_2))$ for any $x_1,x_2\in X$. To get the same 
inequality with $\delta_{\pi_0,U}(r,s)$ replacing $\delta_U(r,s)$, notice that $q\circ r'=q\circ r$ and $q\circ s'=q
\circ s$, for any $r'\in \widetilde{r}$ and $s'\in \widetilde{s}$.
\medbreak

\noindent{\em $q\circ r= q\circ s$.} In this case, both of the summands in the equation (\ref{dist}) will play a role.
The first step is to show that $r$ and $s$ are strictly related:

\begin{lemma}\label{qiso}
Given two descent data $r$ and $s$ relative to a covering $\ucal'=\{U'_i\}_i$ of a complex space $U$, such that $q\circ r= q\circ s$, there is a refinement $\ucal=\{U_i\}_i\succeq\ucal'$ and
collection of holomorphic maps $\phi_i:U_i\to X_1=X\times_\ycal X$ related to $r$ and $s$ by the equations $\partial_0\circ \phi_i=r_i$ and $\partial_1\circ \phi_i=s_i$.
\end{lemma}
\demo
By definition of $Q(\ycal)$, there is a set theoretic map $h_i:U'_i\to X_1$ such that $\partial_0
\circ h_i=r_i$ and $\partial_1\circ h_i=s_i$. By possibly refining $\ucal'$ to $\ucal$, we may suppose that 
\begin{enumerate}
\item[1)] $r_i(U_i)$ is contained in a $\partial_0$ uniformly covered open subspace $V_i\subset X$ for each $i$
and 
\item[2)] each $U_i$ is simply connected.
\end{enumerate}
Let $\partial_0^{-1}(V_i)=\amalg_{k=1}^m W_{ik}$. Then we can partition $U_i$ in $N$ subsets $S_k=h_i^{-1}(W_{ik})$ and there exists an $a$ such that $\overline{S_a}\subset U_i$ has an 
nonempty interior. Define $\phi_i:U_i\to W_{i,a}\subset X_1$ to be the $\partial_0$ lifting of $r_i$. Such a holomorphic map satisfies the equations $\partial_0\circ \phi_i=r_i$ on $U_i$, by construction, and 
$\partial_1 \circ\phi_i=s_i$ on $\overline{S_a}$, hence over all $U_i$, given the property of such a closed subspace and that the equality is between holomorphic maps. 
\enddemo
Notice that the collection $(\phi_i)_i$ does not need to be an isomorphism between the descent data $r$
and $s$, the classifying stacks $\bcal G$ providing an example of this.

\begin{remark}\label{simult}
{\em 
We will need a sharpened version of the lemma just proved: by partitioning $U_i\cap U_j$, if such intersection is nonempty, we may assume that $\overline{S_a}$ is contained in $U_i\cap U_j$. Therefore 
it is possible to lift simultaneously $r_i$ to $\phi_i:U_i\to W_{ia}$ and $r_j$ to $\phi_j:U_j\to W_{ia}$, both $\phi_i$ and $\phi_j$ satisfying the properties of the lemma.
}
\end{remark}

We will assume the covering $\ucal$ associated to the descent data $r$ and $s$ enjoys the properties 
stated in the previous lemma, whose notation we are going to share. There are yet two distinct cases to be considered. Recall that $A\subset X$ is the subspace of $q$ ramification points (see Section \ref{aus}).

\smallbreak

{\em $V_i\cap A=\emptyset$.} We may assume that $\partial_1(W_{ik_1})\cap \partial_1(W_{ik_2})=\emptyset$ if $k_1\neq k_2$ (by possibly further refining $\ucal$). Denote by $M$ the minimum of the distances between pairs of open subspaces $\partial_1(W_{ik})$, $k$ ranging from $1$ to 
$m$. We may assume $M>0$, by shrinking $V_i$ a little if necessary. Now, $r\not\sim s$, so that in particular
$r\neq s$, and we may assume $r_i\neq s_i$ as holomorphic maps $U_i\to X$. Since $q\circ r=q\circ s$  and because of Lemma \ref{qiso}, we conclude that
$$ 
\delta_U(r,s)\geq \sup\limits_{u\in U_i}\bigl\{{\rm d}_X\left(r_i(u),s_i(u)\right)\bigr\}=\sup\limits_{u\in U_i}\bigl\{{\rm d}_X(\partial_1\circ\phi_i (u),\partial_1\circ \phi_i(u))\bigr\}\geq M>0
$$
the same statement can be deduced for $\delta_{\pi_0,U}$ as the same argument applies for any $r'\sim r$
and $s'\sim s$ provided their associated open covering is $\ucal$ or finer.

\smallbreak

{\em $V_i\cap A\neq\emptyset$.} In this case the number $\sup\limits_{u\in U_i}\left\{{\rm d}_X
\left(r_i(u),s_i(u)\right)\right\}$ could possibly be arbitrarly small. However, the only way for $\sup\limits_{u\in 
U_{ij}}\left\{{\rm d}_{X_1}\left(f_{ij}(u),g_{ij}(u)\right)\right\}$ to be arbitrarly small is if $f_{ij}(U_{ij})$ and $
g_{ij}(U_{ij})$ both belong to $W_{ih}$ for some $1\leq h\leq m$, since otherwise that number would be greater or 
equal than $M'=\min\limits_{h_1\neq h_2}\{ {\rm d}_{X_1}(W_{ih_1}, W_{ih_2})\}$ which can be assumed to 
be strictly greater than zero. Let now fix $i,j$ so that $U_i\cap U_j\neq \emptyset$ and consider the worst 
case scenario, i.e. $f_{ij}(U_{ij})$ and $g_{ij}(U_{ij})$ lying in the same $W_{ih}$. We claim that the pair $(\phi_i,
\phi_j)$ forms an isomorphism between $r_{|U_i\cap U_j}$ and $s_{|U_i\cup U_j}$. That being the case, 
since $r\not\sim s$, i.e. $r$ is not isomorphic to $s$, there must be $U_{cd}\neq \emptyset$ such that 
$f_{cd}(U_{cd})\subset W_{ih_1}$ and $g_{cd}(U_{cd})\subset W_{ih_2}$ with $h_1\neq h_2$. This implies that

$$
\delta_U(r,s)\geq \sup\limits_{\stackrel{u\in U_{ij}}{i,j\in\N}}\bigl\{{\rm d}_{X_1}(f_{ij}(u),g_{ij}(u))\bigr\}\geq M'>0
$$
finishing the proof of the proposition.

In order to prove that the pair $(\phi_i,\phi_j)$ is an isomorphism, it suffices to show that $m(f_{ij}, 
\phi_j)=m(\phi_i,g_{ij})$ seen as holomorphic maps defined on $U_i\cap U_j$ and with values in 
$X_1$. This statement follows from knowing that the images of $m(f_{ij},\phi_j)$ and $m(\phi_i,g_{ij})$ lie in the same $W_{ib}$: indeed, this condition on the images implies that $\partial_0$ is injective when restricted 
these images and, since $\partial_0\circ m(f_{ij}, \phi_j)={r_i}_{|U_i\cap U_j}=\partial_0\circ m(\phi_i, g_{ij})$,
we conclude that the  map  $m(f_{ij}, \phi_j)$ must necessarily be the same as $m(\phi_i,g_{ij})$. 
Showing that the images of $m(f_{ij},\phi_j)$ and $m(\phi_i,g_{ij})$ belong to $W_{ih}$ can be further reduced to
the fact (see Remark \ref{simult}) that the ones of $\phi_i$ and $\phi_j$ are both contained in $W_{ia}$. This 
follows by the continuity of $m$: we know that the images of $f_{ij}$ and $g_{ij}$ are in $W_{ih}$ and, by the continuity of $m$ we have that the  maps $z\mapsto m(w_{a},z)$ and $z\mapsto m(z, w_{h})$, for $w_{a}\in W_{ia}$ and $w_{h}\in W_{ih}$, send $W_{ih}$ and 
$W_{ia}$, respectively, to the same $W_{ib}$. 

Consider now two classes $[(r,\phi)]\neq [(s,\psi)]\in \pis_1(\ycal,y)(U)$. These two classes 
being different, in particular it implies that $\widetilde r\neq \widetilde s\in \pis_0(\ycal,y)(U)$ (see Theorem  \ref{pi1}). 
What already proved for $\pis_0(\ycal,y)(U)$, shows that the $\delta_U(r,s)$ summand in the formula 
(\ref{1dist}) is going to be nonzero, hence $\delta_{\pi_1, U}\bigl([(r,\phi)],[(s,\psi)]\bigr)>0$.
\enddemo


\subsection{The function $c(\ycal,-)$}\label{differenz} 
The function $c(\ycal,-)$ associates to every complex space a positive number (possibly $+\infty$). It is defined as follows. Let $U$ be a complex space and 
 $$
 {\sf F}={\sf F}_{\dcal\tms\ucal}=\bigl(F_{ai}:D_a\tms U_i\to \D\tms U , F_{aibj}:D_{ab}\tms U_{ij}\to X_1\bigr)
 $$ 
be a relative analytic disc of $\ycal$ on $U$. We associate a function $\dbb\to \rbb_{\geq 0}$ to ${\sf F}$ as follows.
Take $z\in\dbb$ and a vector $v\in T\dbb$, the holomorphic tangent bundle on $\D$, and consider  $D_a
\times U_i$ with $z\in D_a$. Then, denoting with ${\rm d}F_{ai}$ the differential of the holomorphic map 
$F_{ai}$,  ${\rm d}F_{ai}(z,\xi)v$ is a vector field tangent to $X$ along the points of the image of $(F_{ai})_{z
\times U_i}$,  $\xi$ ranging in $U_i$. Since the  maps $\partial_0$, $\partial_1$ and $m$ are local 
isometries (see Subsection \ref{top}), by differentiating with respect of the variable $z$ the structural equations of ${\sf F}$, $(\star)$ 
and $(\star\star)$ of Subsection \ref{hol-pref}, we notice that the real number 
$$
|{\rm d}F_{ai}(z,\xi)v|:=
q^*H(dF_{ai}(z,\xi)v)
$$
only depends on $\xi, z, v$ and not on the open subspaces $D_a\times U_i$
\footnote{We recall that $H$ denotes the fixed length function on $Q(\ycal)$.}, so $(z,\xi,v)\mapsto
|{\rm d}F_{ai}(z,\xi)v|$ is a well defined real valued continuous function which will be occasionally written as 
$|{\rm d}\fsf(z,\xi)v|$. Moreover, for any $(z,\xi)\in D_{ab}\times U_{ij}$ we have that 

\begin{equation}\label{billcarson}
\vert {\rm d} F_{ai,bj}(z,\xi)v\vert=\vert {\rm d} F_{ai}(z,\xi)v\vert.
\end{equation}
If $K\subset U$ is a compact subset and $z\in\dbb$, we define
 \be\label{FO12}
\vert{\rm d}{\sf F}(z)\vert_K=\sup\limits_{\xi\in K}\sup\limits_{\stackrel{v\in T_z\D}{v\neq 0}}\frac{\vert{\rm d}{\sf F}(z,\xi)v\vert}
{\vert v\vert_{\rm hyp}},  
\ee
where $|v|_{\rm hyp}$ is the length induced by the Poincar\'e metric on $\dbb$, and 
\begin{equation}\label{FO10}
c_{K}(\ycal;\ucal)=\sup\limits_{\stackrel{{\sf F}\in \ycal(\dcal\tms \ucal)}{\mathcal D\in {\rm Cov}\,\D}}\sup
\limits_{z\in\D}\vert{\rm d}{\sf F}(z)\vert_K.  
\end{equation}
Because of the transitivity of the action of Aut$(\dbb)$ and the invariance of the Poincar\'e metric we conclude that
\be\label{FO40}
c_{K}(\ycal;\ucal)=\sup\limits_{\stackrel{{\sf F}\in \ycal(\dcal\tms \ucal)}{\mathcal D\in {\rm Cov}\,\D}}\vert{\rm 
d}{\sf F}(0)\vert_K.
\ee
We finally define
\be\label{FO21}
c\left(\ycal;U\right)=\sup\limits_{\stackrel{K\Sbs U}{\ucal\in {\rm Cov}\,U}}c_{K}\left(\ycal;\ucal\right).
\ee
Since $\partial_0$ and $\partial_1$ are local isometries, if the relative analytic disc 
\be\label{disc2}
{\sf G}={\sf G}_{\dcal\tms\ucal}=\bigl(G_{ai}:D_a\tms U_i\to \D\tms\cbb , G_{aibj}:D_{ab}\tms 
U_{ij}\to X_1\big)
\ee
is equivalent to ${\sf F}$, i.e. they have the same class in $\pis_0(\ycal,y)(\dbb\times U)$, the functions $|{\rm d}
\fsf(z,\xi)v|$ and $|{\rm d}\gsf(z,\xi)v|$ coincide. 

Given a pair $[\fsf,\Phi]$ representing a class of $\pis_1(\ycal,y)
(\dbb\times U)$, where $\fsf$ is a relative analytic disc on $U$, and $\Phi=\{\Phi_{ai}\}_{ai}$ is an automorphism of $\fsf$ (see Section \ref{hol-pref}), for each compact set $K\subset U$ we have well defined functions $|{\rm d}\fsf(z)|_K$ and
$|{\rm d}\Phi(z)|_K$, as in equation (\ref{FO12}). Keeping in mind the Theorem \ref{pi1}, we deduce that 
$|{\rm d}\fsf(z)|_K=|{\rm d}\Phi(z)|_K$. If $[\gsf,\Psi]$ is a pair equivalent to $[\fsf,\Phi]$, that is their images coincide 
in $\pis_1(\ycal,y)(\dbb\times U)$, then
\be\label{disc34}
\vert{\rm d}\,{\sf F}(z)\vert_K=\vert{\rm d}\,{\Phi}(z)\vert_K=\vert{\rm d}\,\Psi(z)\vert_K=\vert{\rm d}\,{\sf G}(z)
\vert_K.
\ee

\begin{remark}\label{costante}
{\em Because of the last observations, the vanishing of the $H$-norm of the derivative of a descent data or 
of one of their isomorphisms is equivalent to their relevant classes in $\pis_0$ or $\pis_1$ being constant (see Definition \ref{sez-cost}).}
\end{remark}

Under the same notation as in Subsection \ref{top}, we prove the following fundamental
\begin{lemma}\label{lcat1}
Let $\alpha_1,\alpha_2\in \pis_i(\ycal,y)(U)$, for $i=0,1$ and $\csf_{\alpha_1,\alpha_2)}$ be  an analytic chain through $\alpha_1$ and $\alpha_2$. Then 
\begin{equation}\label{EQ33}
{\delta}_{\pi_i,U}(\a_1,\a_2)\le 2\, c(\ycal;U)\, l({\sf C}_{(\a_1,\a_2)})
\end{equation}
(see the equation (\ref{lcat}) for the definition of length of an analytic chain). In particular,
\be\label{lcat45}
{\rm d}_{\rm Kob}(\a_1,\a_2)\ge\frac{{\delta}_{\pi_i,U}(\a_1,\a_2)}{2\, c(\ycal;U)}.
\ee
for $i=0,1$ where the left end side has been defined in the equation (\ref{ddisc}) with $\pcal=\ycal$.
\end{lemma}

\demo
It is sufficient to prove these statements for analytic discs instead of chains. Let $\alpha_1\neq\alpha_2\in \pis_0(\ycal,y)(U)$ and $r=({r}_i,
{f}_{ij})\neq s=({s}_i,{g}_{ij})$ be two descent data in $\a_1$, $\a_2$ respectively. Let $\fsf=\bigl( F_{ai},F_{aibj}\bigr)$ be an analytic disc such that $r=z_1^*\fsf$ and $s=z_2^*\fsf$ for $z_1$ and $z_2$
points which we initially suppose are in a $D_a$ open of a covering of  $\dbb$ associated to $\fsf$. 
Using the notation introduced in Subsection \ref{dischi}, in particular we have that 
$F_{ai}(z_1,u)=r_{ai}(u)$ and $F_{ai}(z_2,u)=s_{ai}(u)$. Let $\gamma:[0,1]\to D_a$ a geodesic arc
in $\dbb$, endowed by the Poincar\'e metric, such that $\gamma(0)=z_1$, $\gamma(1)=z_2$ and, for  
$u\in U$ fixed, we call $\gamma_u$ the path $t\to F_{ai}(\gamma(t),u)$.  For each $u\in U_i$ we have
\begin{eqnarray}\label{lun1} 
{\rm d}_X(r_{ai}(u),s_{ai}(u))\le {\rm length}(\gamma_u)\,\!\!\!\!&=&\!\!\!\!\int_{0}^{1}H\left(\,{\rm d} F_{ai}(\gamma(t),u)\cdot\gamma'(t)\, \right){\rm d}t\\
&\le& c(\ycal;U)\int_{0}^{1}\vert\gamma'(t)\vert_{{\rm hyp}}{\rm d}t
= c(\ycal;U)\r_{{_\D}}(z_1,z_2)\nonumber.
\end{eqnarray}
The same argument applies to the transition functions $f_{aibj}$ and $g_{aibj}$ and, using the remark in
equation (\ref{billcarson}), we get
\begin{eqnarray}\label{lun2} 
{\rm d}_{X_1}(f_{aibj}(u),g_{aibj}(u))&\le&{\rm length}(\,\g_u)=
\int_0^1H\left({\rm d} F_{aibj}(\gamma(t),u)(\gamma'(t))\right){\rm d}t\le\\\nonumber
&& c(\ycal;U)\int_{0}^{1}\vert\gamma'(t)\vert_{{\rm hyp}}{\rm d}t= c(\ycal;U)\r_{{_\D}}(z_1,z_2)\nonumber.
\end{eqnarray}
Taking the supremum with respect to $u$ and the indices $i, ai$, we deduce from (\ref{lun1}) and 
(\ref{lun2}) that 
$$
\delta_U(r,s)\le 2 c(\ycal;U)\r_{{_\D}}(z_1,z_2)
$$
which implies 
$$
{\delta}_{\pi_0,U}(\a_1,\a_2)=\inf\limits_{\stackrel{r\in\a}{s\in\b}}\delta_U(r,s)
\le 2 c(\ycal;U)\r_{{_\D}}(z_1,z_2).
$$
In the general case, when $z_1$ and $z_2$ do not belong to the same open $D_a$, we partition  
 $\gamma([0,1])\subset\dbb$ with points $z_1=\z_1=\gamma(t_1), \z_2=\gamma(t_2), \cdots , \z_m=z_2
\gamma(1)$ so that  the geodetic arcs $\gamma(\z_l,z_{l+1})$ are contained in $D_{a_0},D_{a_1},\cdots
D_{a_{m-1}}$, respectively, for $0\leq l\leq m-1$. Then we argue as before on each arc. This proves 
th equation (\ref{EQ33}).

The statement regarding $\pis_1(\ycal,y)(U)$ is proved in the same fashion and by using  the equalities
in $(\ref{disc34})$. 
\enddemo


\section{Kobayashi hyperbolicity implies Brody hyperbolicity}\label{kob-sec}
We can now proceed with the proof of the Brody theorem for stacks. Classically, that theorem refers
to the implication ``compactness and Brody hyperbolicity imply Kobayashi hyperbolicity'', the other 
implication holding in general and simply a consequence of non Kobayashi hyperbolicity of $\cbb$ and
that every holomorphic map is a contraction with respect to the Kobayashi presudodistance.
In the context of stacks, even this simple implication, though not using much of the theory developed
in the last sections, is not entirely obvious and this section is devoted to its proof. 

\begin{theorem}\label{KB}
Let $\ycal$ be a Kobayashi hyperbolic analytic stack. Then $\ycal$ is Brody hyperbolic. 
\end{theorem}
\demo 
Suppose that $\pis_0(\ycal,y)$ is not Brody hyperbolic; then there exists a section $r\in \pis_0(\ycal,y)(\cbb\times U)$
not in the image of $p^*: \pis_0(\ycal,y)(U)\to \pis_0(\ycal,y)(\cbb\times U)$, $p$ being the projection.
Let us introduce the notation $r=(r_i(z,u), f_{ij}(z,u))$ over a covering $\{V\}_{i}$ 
of $\cbb\times U$; notice that we may take $V_i$ to be of the kind $D_{k(i)}\times U_{\alpha(i)}$ for a disc
$D_{k(i)}$ centered at some complex number $z\in\cbb$ and open set $U_{\alpha(i)}$ of $U$.
 We wish to construct two sections, or objects, $r_1\not\cong r_2\in$ Ob$(\ccal\xdot(V))$ for some complex
space $V$, whose Kobayashi pseudodistance is zero. Take $V=\cbb\times U$ and consider the  
two sections: $r_1=s$ and $r_2=p^*(i_0^*(s))$, where $i_0:U\to\cbb\times U$ is the 
embedding in zero. By assumption, $r_1\not\cong r_2$. To show that the Kobayashi pseudodistance between
$r_1$ and $r_2$ is zero, we construct relative analytic chains (see Subsection \ref{dischi}) $s_n^\lambda$ between
them for $n\in\nbb$, where $\lambda\in \dbb_2=\{z\in\cbb: |z|<2\}$. For each $n$, $s^\lambda_n$ is, in fact, a relative 
analytic disc. We define first the descent data $s_1^\lambda$ on $\dbb_2\times \cbb\times U$; to describe the open 
covering we use its cross sections for $\lambda=$const as $t_{1/\lambda}(D_{k(i)})\times U_{\alpha(i)}$, where
$t_\omega:\cbb\to \cbb$ is the automorphism $z\to\omega z$. For $\lambda=0$ the open covering of $\cbb\times 
U$ is $\{\cbb\times U_{\alpha_0}\}_{\alpha_0}$, where $\wcal=\{U_{\alpha_0}\}_{\alpha_0}$ is the covering of $U$ given by the open sets appearing 
in the covering of $\{0\}\times U$ induced by the $V_i$. 
The holomorphic map is $(\lambda, z, u)\mapsto r_{i}(t_{1/\lambda}^{-1}(z),u)=(\lambda z, u)$ for $z\in t_{1/\lambda}
(D_{k(i)})$ and $u\in U_{\alpha(i)}$. The transition functions are similarly defined as $(\lambda, z,u)\mapsto f_{ij}(\lambda 
z, u)$. Notice that $s^1=r_1$ and $s^0=r_2$.
The general relative analytic disc $s_n^\lambda$ is associated to the open covering having open sets
$\{t_{\frac{1}{n\lambda}}(D_{k(i)})\times U_{\alpha(i)}: \lambda\in \dbb_2\}$ and holomorphic maps $r_i(n\lambda z,u)$
and $f_{ij}(n\lambda z,u)$ for $\lambda\in \dbb_2$, $z\in t_{\frac{1}{n\lambda}}(D_{k(i)})$ and $u\in U_{\alpha(i)}$.
For each $n$, we have $s_n^{1/n}=r_1$ and $s_n^0=r_2$, thus (see equation (\ref{ddisc})
$$
{\rm d}_{\rm Kob}(r_1, r_2)\leq \rho_{\dbb_2}(1/n, 0)=\frac{1}{2}\log\Big (\frac{2n+1}{2n-1}\Big)
$$  
which tends to zero as $n$ goes to infinity.

Consider now the case of $\pis_1(\ycal,y)$ not being Brody hyperbolic, thus there exists a section $\phi\in \pis_1(\ycal,y)(\cbb\times U)$ not in $p^*(\pis_1(\ycal,y)(U))$.  We recall (see Subection \ref{hol-pref}) that, 
for a complex space $V$, the sections in $\pis_1(\ycal,y)(V)$ are represented by classes $[(\psi_i)_i]$
of automorphisms of some descent datum $(r_i,f_{ij})$. In particular, $\psi_i:V_i\to X\times_\ycal X$ are holomorphic maps from open sets $V_i$ of an open covering $\vcal$ of $V$, satisying certain identities.
From the section $\phi$ just mentioned, we produce two different sections $\theta_1$ and $\theta_2$ in $
\pis_1(\ycal,y)(\cbb\times U)$ as in the $\pis_0$ case: $\theta_1=\phi$ and $\theta_2=p^*(i_0^*(\phi))$.
Notice that $\theta_2$ is an automorphism of a descent datum $r_2\in${\sf Ob}$(\ccal\xdot(\cbb\times U))$ which
lies in $p^*(${\sf Ob}$(\ccal\xdot(U)))$ and $\theta_1$ is an automorphism of a descent data $r_1$
which may possibly be in the image of $p^*$, although $\theta_1$ itself does not. However, it necessarily has to be $r_1\neq r_2\in \pis_0$ given the hypothesis on $\phi$. Keeping this in mind, to denote the descent data 
$r_1$ and $r_2$ we will use the notation already 
introduced  during the discussion of the $\pis_0$ case: $r_1=(r_i(z,u), f_{ij}(z,u))$ and $r_2=(r_i'(z,u),f_{ij}'(z,u))$, $r'_i(z,u)=r(0,u)$ and $f_{ij}'(z,u)=f_{ij}(0,u)$.  We wish to find an  
analytic disc $\Omega$ relative to $\cbb\times U$ between $\theta_1$ and $\theta_2$. 
Define $\Omega: \dbb_2\times (\cbb\times U)\to \pis_1(\ycal,y)$ by 
$$
\Omega(\delta_U,z,u)=\{(\Omega_k)_k\; \text{such that}\; \Omega_k(\delta,z, u_k)=\phi_k(\delta z, u_k)\},
$$
where $\phi=(\phi_k)_k$ and $\phi_k:W_k\to X_1$ are holomorphic maps defining the extension of $\phi$ over the open covering $\wcal$ of $U$, 
precedently introduced. For any fixed $\delta\in \dbb_2$, $m(\phi_h, f^\delta_{hk})(z, u)=m(f^\delta_{hk}, 
\phi_k)(z,u)$, which shows that $\Omega_{|\{\delta\}\times\cbb\times U}$ is an automorphism of $s^\delta$ for any 
$\delta\in \dbb_2$ and $\Omega$ indeed induces an analytic disc in $\pis_1(\ycal,y)(\cbb\times U)$. Also, 
$\Omega_{|\{0\}\times\cbb\times U}=\theta_1$ and $\Omega_{|\{1\}\times\cbb\times U}=\theta_2$ so that it is an 
analytic disc between $\theta_1$ and $\theta_2$. The sequence of analytic discs $\Omega^n$ relative to $\cbb
\times U$ is now introduced in analogy with the $\pis_0$ earlier discussion and the vanishing of the Kobayashi pseudodistance between $\theta_1$ and $\theta_2$ follows at once.
\enddemo

\section{Compactness and Brody hyperbolicity imply Kobayashi hyperbolicity}\label{dura}

Brody and Kobayashi hyperbolicity of an analytic stack $p:X\to \ycal$ are statements concerning the holotopy presheaves $\pis_0(\ycal,y)$ and $\pis_1(\ycal,y)$, as conceived in definitions \ref{br-campi} and 
\ref{kob-iper}. Surprisingly, it turns out that for compacts  Deligne-Mumford analytic stacks, the Brody hyperbolicity content 
in the holotopy presheaf $\pis_0$ absorbed the one of $\pis_1$, to the point of making the latter irrelevant.

\begin{theorem}\label{BRO}
Let $\ycal$ be a compact,  Deligne-Mumford analytic stack. If $\pis_0(\ycal,y)$ is Brody hyperbolic, then
\begin{itemize}
\item[i)] for any complex space $U$ we have $c(\ycal,U)<+\infty$;
\item[ii)] $\pis_i(\ycal,y)$ are Kobayashi hyperbolic, for $i=0,1$. 
\end{itemize}
\end{theorem}

\demo
The statement i) implies statement ii). Indeed, let $\alpha_1\neq \alpha_2\in \pis_i(\ycal,y)(U)$ admissible
sections. From i) and Lemma \ref{lcat1}, we have
$$
{\rm d}_{{\rm Kob}}(\a_1,\a_2)\ge\frac{\delta_{\pi_i,U}(\a_1,\a_2)}{2c(\ycal;U)}>0
$$ 
hence $\pis_i(\ycal,y)(U)$ are Kobayashi hyperbolic.
The rest of the section will be devoted to the proof of the first assertion, which will be split in few steps.   
\smallbreak
Assume by contraddiction that $c(\ycal, U)=+\infty$ for some complex space $U$. Then 
there exists a sequence $\{\ucal_\nu\}_\nu$ of coverings of $U$ associated to a sequence of 
relative analytic discs $\{\fsf^\nu\}_\nu$, i.e. descent data, over $\dbb\times U$ and points $\{u_\nu\}_\nu\sbs U$
all with the property that $\lim_{\nu\to +\infty}|{\rm d}\fsf^\nu(0,u_\nu)|= +\infty$ (see equation (\ref{FO40})). 
The descent datum $\fsf^\nu$ induces on $\dbb\times \{u_\nu\}$, hence on $\dbb$, consequently $\{\fsf^\nu\}_\nu$ induces on $\D$ a sequence of descent data which we will keep writing in the same way with the property that $\lim_{\nu\to +\infty}|{\rm d}\fsf^\nu |=+\infty$. This reduces the argument
to $U=\{{\sf pt}\}$. 

\begin{lemma}\label{fun}
If $c(\ycal,U)=+\infty$, then there is a sequence $\{f^\nu\}_\nu$ of holomorphic maps $f^\nu:\dbb\to 
Q(\ycal)$ (see Section \ref{aus}) such that $\lim_{\nu\to +\infty}|\drm f^\nu(0)|=+\infty$. 
\end{lemma}
\demo
There exists a commutative diagram
\be\label{fatt}
\xymatrix{\ccal\xdot(U)\ar[rr]\ar[rd]_{\phi}& &\pi_0(\ycal)(U)\ar[ld]^{\phi_Q}\\ &{\rm Hol}(U,Q(\ycal))&}
\ee
where the application $\phi$ is defined by associating to a descent datum $r=(r_i, f_{ij})$ on $U=\cup_{i\in\N}U_i$ the 
holomorphic map $f_r:U\to Q(\ycal)$ defined for $u\in U_i$ as $f_r(u)=q(r_i(u))$, where $q:X\to Q(\ycal)$ 
denotes as usual the holomorphic projection. The map $f^r$ is well defined since if $u\in U_i\cap U_j$, then $q(r_i(u))=
q(r_j(u))$, since there exists $w=f_{ij}(u)\in X\times_\ycal X$ such that $\partial_0(w)=r_i(u)$ and $\partial_1
(w)=r_j(u)$. If two descent data $r$ and $s$ represent the same class in $\pis_0(\ycal,y)(U)$, then $\phi(r)=
\phi(s)$, thus there is a well defined application $\phi_Q$ making the diagram (\ref{fatt}) commutative. 
In particular, to each relative analytic disc $\fsf$ on $\dbb\times U$ we have a holomorphic map 
$f:=\phi(\fsf):\dbb\times U\to Q(\ycal)$ and a real continuous function $|\drm f(z)|$.  From the very definition
we have that $|\drm\fsf(z)|=|\drm f(z)|$, thus the statement of the lemma follows at once by 
taking $f^\nu=\phi(\fsf^\nu)$.  
\enddemo

The sequence of the maps
$f^\nu:\dbb\to Q(\ycal)$ may be reparametrized, by means of the ``reparametrization Lemma'' (cfr. \cite{brody}), to
get a sequence of maps $\widetilde f^\nu:\dbb_\nu\to Q(\ycal)$, where $\dbb_\nu=\{|z|<\nu\}$ and $|{\rm d}\widetilde f^\nu(0)|=1$ for all $\nu$. By Ascoli-Arzel\acc Theorem, there exists a subsequence $\widetilde f^\mu$ uniformly convergent on compacts to a holomorphic map $f:\cbb\to Q(\ycal)$, which is not constant since $|{\rm d}f(0)|=1$ (cfr. \cite{brody}).  
In the rest of the section, the reparametrized maps $\widetilde f^\nu$ will  be denoted $f^\nu$ and the domain $\dbb_\nu$ will be specified to distinguish them from those defined on $\dbb$. 

Recall (Section \ref{aus}) that $A$ is the ramification locus of $q:X\to Q(\ycal)$ and $B=q(A)\subset Q(\ycal)$ 
is the branch locus.  One of the following three cases necessarily occurs

\begin{itemize}
\item[Case 1)] $f(\C)\cap B=\emptyset$;
\item[Case 2)] $f^{-1}(B)$ is a discrete set of points of $\cbb$;
\item[Case 3)] $f(\C)\sbs B$.
\end{itemize}
We would like to ``invert'' the map $\phi$, i.e. find an application $\psi$ providing a descent datum
in $X$ associated to a holomorphic map $U\to Q(\ycal)$ such that $\phi\circ \psi={\sf id}$.  This property
implies that if $f$ is nonconstant, $\psi(f)$ is a nonconstant descent datum on $\cbb$ (see Definition \ref{sez-cost}) and $[\psi(f)]\in \pis_0(\ycal,y)(\cbb)$ is a nonconstant class (see Remark \ref{costante}).

We are going to consider the existsence of $\psi$ in the above three cases separately, starting with the first.

\subsection{Case 1}

\begin{prop}\label{soll}
Let $h:\cbb\to Q(\ycal)\ssmi B$ be a holomorphic map. 
Then there is a descent datum $\psi(h)$ in $X$ associated to $h$, such that $\phi(\psi(h))=h$. 
\end{prop}

\demo
Let $U$ be any complex space; we will take $U=\cbb$ just at the end of the proof.
Cover $U$ with simply 
connected open subspaces $U_i$ and lift the restrictions of $h$ to continuous maps $h_i':U_i\to X$. By 
\cite{grezzo}, we know that $q:X\ssmi A\to Q(\ycal)\ssmi B$ is a local biholomorphism, hence the $h_i'$
are holomorphic maps. Let $U_i\cap U_j\neq \emptyset$; to relate $h_i'$ and $h_j'$ on $U_{ij}=U_i\cap U_j$
we argue as follows. Since $h_{|U_{ij}}=h_{|U_{ij}}$, there exists a set theoretic map $\omega: U_{ij}
\to X_1$ such that $\partial_0\circ \omega={h_i'}_{|U_{ij}}$ and $\partial_1\circ \omega={h'_j}_{|U_{ij}}$. Refining
the covering $\{U_{i}\}_i$, if necessary, we can assume the $U_i$ are relatively compact, $U_{ij}$ are simply connected and $h_i'(U_{ij})$ are uniformly covered by $\partial_0$ so that $\partial_0^{-1}h_i'(U_{ij})=
\amalg_{k=1}^l V_k$ and $\partial_0:V_k\to U_{ij}$ is a biholomorphism for any $k$. Consider first the case of 
a one dimensional complex space $U$. As the codomain of $\omega$ has a finite number of connected components, 
there is an index $a\leq n$ such that $\omega^{-1}(V_a)=L$ contains uncountably many points, hence there is an 
accumulation point among them. Call $f_{ij}':U_{ij}\to V_a$  the $\partial_0$ lifting of ${h_i'}_{|U_{ij}}$ to $V_a$; it 
coincides with $\omega$ on $L$ since ${\partial_0}_{|V_a}$ is injective. Thus, $\partial_0\circ f'_{ij}={h'_i}_{|U_{ij}}$ 
on $U_{ij}$, $\partial_1\circ f'_{ij}={h'_j}_{|U_{ij}}$ on $L$,  hence on all $U_{ij}$. This shows that $f'_{ij}$ is a 
transition function.  

Let us assume inductively to be able to define the transition functions on each complex space $U$ of dimension 
$n-1$. We think of $U_{ij}$ to be embedded in $\cbb^N$ and consider the intersections
$U_{ij}\cap H_\lambda$, $H_\lambda$ ranging among all $(N-1)$ dimensional hyperplanes  in $\cbb^N$.  
On each of these intersections we have transition functions with values in $V_a$ where the index $a$ may depend 
on $\lambda$, but being the possible indices finite, there are infinitely many $\mu$ such that $U_{ij}\cap H_
\mu\neq \emptyset$ and whose transition functions have values in $V_a$ for a fixed $a$. Call $l_{ij}$ the 
holomorphic map on $L:=\bigcup_\mu U_{ij}\cap H_\mu\to V_a$ defined by taking on each $U_{ij}\cap H_\mu$
the transition function on it; notice that $l_{ij}$ is well defined since all the transition functions are $\partial_0$ liftings to $V_a$ and ${\partial_0}_{|V_a}$ is injective. Call now $f'_{ij}:U_{ij}\to V_a$ the $\partial_0$ lifting of 
${h'_i}_{|U_{ij}}$ to $V_a$. As expected, $\partial_0\circ f'_{ij}={h'_i}_{|U_{ij}}$ and $\partial_1\circ
f'_{ij}={h'_j}_{|U_{ij}}$ coincide, as functions defined on $U_{ij}$, on an analytic subspace containing $L$, but the only such a space is $U_{ij}$ itself.
A good candidate for the descent datum $\psi(h)$ could be $(h_i', f_{ij}')$, which indeed looks very much like a 
descent datum on $U$, except  that, in general, the $f'_{ij}$ do not satisfy the cocyle condition 
$m(f_{ij}',f_{jk}')=f_{ik}'$ (see Definition \ref{discesa}), $m$ being the multiplication of $\xdot$.  The idea is the one
of trying to replace $f'_{ik}$ with $m(f_{ij}', f_{jk}')$, the problem being that the latter is defined only on the 
triple intersection $U_{ijk}$. 

\begin{lemma}
Given a non constant holomorphic map $h:U\to Q(\ycal)$, refer to $h'_i:U_i\to X_1$ and $f'_{ij}:U_{ij}\to
X_1$ as the holomorphic maps obtained by means of the above constructions. Then there exists a refinement
$\{B_\alpha\}_\alpha$ of $\{U_i\}_i$ and transition functions $\{g_{\alpha\beta}\}_{\alpha\beta}$ such that, if $B_i$, 
$i=1,2,3$ are open sets with nonempty triple intersection, for any chosen transition functions $g_{12}$ and $g_{23}$, 
there exists a holomorphic map $g_{13}:B_{13}\to X_1$ satisfying the cocycle condition $m(g_{12}, 
g_{23})=g_{13}$ on $B_{123}$. Moreover, $\phi((h'_\alpha, g_{\alpha\beta}))=h$.

\end{lemma}

\demo 
Denote the holomorphic map $m(f'_{12}, f'_{23}):U_{123}\to X_1$ by $l_{13}$, choose points $u\in U_{123}$ 
and $p\in U_{13}$ and a path $\gamma: [0,1]\to U_{13}$ connecting $u$ to $p$.  We would be done if $h_1'\circ 
\gamma$ lifted to a $\tilde{\gamma}$ starting in $l_{13}(u)\in X_1^{13}$ where $X_1^{13}$ is the connected component of 
$X_1$ containing the image of $l_{13}$, but this may necessarily not be the case since $h_1'(U_{13})$ might not be 
contained in $\partial_0(W_1^{13})$. The possibility of lifting such a map only depends on $p$, $X$ and $\ycal$, 
hence, using that $\partial_0:X_1\to X$ is a covering, we can refine $\{U_i\}_i$ to a covering $\{B_\alpha\}_{\alpha\in\nbb}$ of $U$ with the following properties: 
\begin{enumerate}
\item[1)] $\{B_\alpha\}_\alpha$ is a locally finite refinement of $\{U_i\}_i$;
\item[2)] the intersections $B_{\alpha\beta}$ are simply connected;
\item[3)] if $\tau(\alpha)$ is less than the integers $i$ such that $B_\alpha\subset U_i$, for each triple intersection $B_{\alpha\beta\gamma}\neq \emptyset$ with $\alpha<\beta<\gamma$, we have $h'_{\tau(\alpha)}(B_{\alpha\beta})$, $h'_{\tau(\alpha)}(B_{\beta\gamma})$ and $h'_{\tau(\alpha)}(B_{\alpha\gamma})$ are contained in the same $\partial_0$ uniformly covered open subspace of $X$.
\end{enumerate}
We now create a descent datum $r$ on $U$ using the covering $\{B_\alpha\}_\alpha$ and the holomorphic maps
$h'_i, f_{ij}'$. For each $\alpha$ define $r_\alpha={h'_{\tau(\alpha)}}_{|B_\alpha}$ and let the transition functions $g_{\alpha\beta}:B_{\alpha\beta}\to X_1$ be 
$$
g_{\a\b}=\begin{cases}{\rm id},& {\rm se}\>B_\a\cup B_\b\sbs U_{\tau(\a)}\\ {f'_{\tau(\a)\tau(\b)}}_{|_{B_\a\cap B_\b}}, &{\rm otherwise}.
\end{cases}
$$ 
The third condition on the covering $\{B_\alpha\}_\alpha$ implies that, if $B_{\alpha\beta\gamma}$ is non empty,
${h'_\alpha}_{|B_{\alpha\beta}}$ lifts to a holomorphic map $\tilde{h}_\alpha$ with image in the same connected 
component of $X_1$  as the one of the image of $m(g_{\alpha\beta}, g_{\beta\gamma})$. Since $\partial_0\circ 
\tilde{h}_\alpha=\partial_0\circ m(g_{\alpha\beta}, g_{\beta\gamma})$, by the injectivity of $\partial_0$ on that 
connected component, we conclude that $\tilde{h}_\alpha$ is an extension of $m(g_{\alpha\beta}, g_{\beta
\gamma})$ to all $B_{\alpha\gamma}$ and it can be taken to be the transition function $g_{\alpha\gamma}$ satisfying the cocyle condition.
\enddemo
A recursive argument can extend this construction to the quadruple and higher intersections,  the only possible 
indeterminacy being in how to choose the transition functions: for instance, $g_{14}$ could be the extension
of $m(g_{12}, g_{24})$ or $m(g_{13}, g_{34})$. However, using the cocycle relations 
\bit
\item[1)] $m(g_{12},g_{23})=g_{13}$, 
\item[2)] $m(g_{23},g_{34})=g_{24}$  
\eit
and assuming to have defined $g_{14}$ as the extension of $m(g_{12},g_{24})$, we check that it verifies the relation 
 $g_{14}=m(g_{13}, g_{34})$, as well: from the first equation the left end side of the equation is $m(m(g_{12},
 g_{23}), g_{34})$ and from the second relation we have that this is precisely $m(g_{12,}g_{24})$, which by definition is precisely $g_{14}$. 
\hfill $\square$

To finish the proof of the proposition \ref{soll} we notice that we can cover $\cbb$ with open subspaces $\{B_i\}_i$
in such a way $B_j\cap \bigcup_{i<j}B_i$ is contractible, therefore there cannot be conflicting
conditions on the function $g_{13}$ of the previous lemma for such a covering. \hfill $\square$

\smallbreak

This finishes the construction of the descent data $\psi(h)$ associated to a holomorphic map $h: U\to Q(\ycal)$, 
hence of the application $\psi$ and the proof of Proposition \ref{soll} and Case 1 are settled. 
\begin{remark}\label{noncostante}
{\em Notice that if the holomorphic map $h$ is nonconstant, then all the liftings $h_i'$ are nonconstant. Using results proved in Subsection \ref{differenz}, we see that we can check whether a section in $[(h_i', h_{ij}')]\in \pis_0(\ycal,y)(\cbb)$ is nonconstant or not, simply by showing that the real valued function $z\mapsto |{\rm d} h_i'(z)|$ is nonzero for $z\in\cbb$. It follows that, if $h$ is nonconstant, so does  $[(h_i', h_{ij}')]$. }
\end{remark}

\smallbreak


\subsection{Case 2} 
Without loss of generality, we can assume that each connected component of the atlas $Z=\amalg_{i=1}^NZ^{(i)}$, containing $X$ as a relatively compact open subspace, is Stein and that $X^{(i)}=\{\zeta\in Z^{(i)}: \phi_i(\zeta)< c\}$,
where $\phi_i:Z^{(i)}\to \rbb$ is a strictly plurisubharmonic exaustion function for $Z^{(i)}$. 
\smallbreak

By assumption $f^{-1}(B)$ is a discrete subset $S$ of $C$. 
We will first show that for each $\zeta\in S$ there is a disc $\Delta_\zeta$ centered in $\zeta$ having the  property that $f_{|\Delta_\zeta^*}$ lifts to a holomorphic map $\tilde{f}_\zeta:\Delta_\zeta^*\to X$, where $\Delta_\z^*$ is the 
punctured disc $\Delta_\zeta\ssmi\{\zeta\}$. This follows if we prove that we can lift $f_{|\Delta_\z^*}$ to a continuous map in $X$ which is equivalent to showing that 
\be\label{soll2}
f_\ast\left(\pi_1(\Delta^\ast_\zeta,z)\right)\sbs q_\ast\, \pi_1\left(X\ssmi A,x)\right)\sbs\pi_1\left(Q(\ycal)\ssmi B,q(x)\right)
\ee
where $q(x)=f(z)$. Let $\gamma$ be a generator of $\pi_1(\cbb\ssmi S,z)$, i.e. a simple closed path around one of the points of $\zeta\in S$. Since the holomorphic maps $f^\nu=\psi(\fsf^\nu)$ converge to $f$ 
uniformly on compact sets and $Q(\ycal)$ is locally topologically contractible, there exists an index $\mu$ such that  
\begin{enumerate}
\item[1)]\label{prim} $f^\mu\circ\gamma$ is homotopic to $f\circ \gamma$ and $(f^\mu\circ\gamma)([0,1])\cap B=\emptyset$;
\item[2)] $\gamma([0,1])\subset \dbb_\mu$.
\end{enumerate}
The map $f^\mu$ can be lifted to $X$ if restricted to an appropriate punctured open disc $\Delta_\z^*$,  $\Delta_\z$ being a contractible open neighborhood of $\z$ in $\cbb$: just take the holomorphic map $F^\mu_i$ part of the descent datum $\fsf^\mu=(F^\mu_i, F^\mu_{ij})$, which is defined on $\gamma([0,1])$. We can assume that $\gamma$ is close enough to $\zeta$ to have image contained on some open set $U_i\ni\zeta$ associated to $\fsf^\mu$. Let $\Delta_\zeta$ be such that $\gamma([0,1])\subset\Delta_\zeta\subset U_i$. Now consider the closed path $F^\nu_i\circ \gamma$ in $X$: by construction and condition \ref{prim} on $\mu$ in particular it induces a homotopy class $[F^\nu_i\circ\gamma]\in \pi_1(X\ssmi A,x)$; since $q\circ(F^\nu_i\circ \gamma)=f^\nu\circ \gamma$ we have

$$
[f\circ\gamma]=[f^\mu\circ\gamma]=q_*(\overline{F^\mu_i\circ\gamma})\in\pi_1(Q(\ycal),y)
$$
which means that $f_*([\gamma])$ lies in $q_*\pi_1(X\ssmi A,x)\subset \pi_1(Q(\ycal)\ssmi B, f(z))$, and so 
does \\ $f_*\pi_1(\cbb\ssmi S,x)$. We can then lift $f_{\Delta^*_\zeta}$ to a holomorphic map $\widetilde{f}_\zeta:
\Delta^*_\zeta\to X\ssmi A$, since $q$ is unramified on $X\ssmi A$. 

Take $W$ to be an open subspace of $\cbb\ssmi S$ such that $\{W, \Delta^*_\zeta: \zeta\in S\}$ is an 
open covering of $U$.
Since $f(W)\subset Q(\ycal)\ssmi B$, by the Case 1 of the Theorem \ref{soll}, $f_{|W}$ lifts to a descent datum 
$r=\psi(f_{|W})$. The datum $r$ is extended to a descent datum on the whole $\cbb\ssmi S$ by using the liftings $\tilde{f}_
\zeta$ as holomorphic maps on $\Delta^*_\zeta$ and creating transition functions as shown in the construction of 
the application $\psi$ in the unramified case. To further extend $r$ to all of $\cbb$, we just have to extend the domain
of the  homolorphic maps $\widetilde{f}_\zeta$ from $\Delta^*_\zeta$ to the whole $\Delta_\zeta$. Let $X^{(1)}\Sbs Z^{(1)}$ be the connected components of the atlases $X$ and $Z$ containing the image of $\widetilde{f}_\zeta$.
Embedding $Z^{(1)}$ as a closed analytic subspace in $\cbb^m$, the map $\widetilde{f}_\zeta$ is expressible in terms
of $m$ complex valued holomorphic functions $f_1,f_2,\cdots, f_m$ which are bounded given the compactness of $\oli X^{(1)}$, thus each of them extendable on $\Delta_\zeta$ by Riemann Theorem. Call $\widehat{f}_\zeta$ the extension of 
$\widetilde{f}_\zeta$ to $\Delta_\zeta$. Suppose now that $\widehat{f}_\zeta(\zeta)$ is a point of the boundary of $X^{(1)}$
in $Z^{(1)}$, then the nonconstant plurisubharmonic function $\phi_1\circ \hat{f}_\zeta$, where $\phi_1$ is the exaustion function introduced before, has a maximum in $\zeta$,
which is absurd. We conclude that the image of $\widehat{f}_\zeta$ is contained in $X^{(1)}$ and the collection of the functions $\widehat{f}_\zeta$, along with the rest of the data in $s$, yields the sought extension of $s$ to $\cbb$. Because $|{\rm d}f(0)|=1$, $f$ is nonconstant, therefore as pointed in Remark \ref{noncostante}, the class of this descent datum in $\pis_0(\ycal,y)(\cbb)$ is nonconstant, as well.  

\smallbreak

\subsection{Case 3}
We assume now that $f(\cbb)\subset B$. Moreover, since $q$ is quasi finite and $B$ is analytic (\cite{grezzo}), there exists a closed analytic subset $B_1$ of $B$ such that $\dim_x(B_1)<\dim_x(B)$ for all $x\in B$ and $q$ is a local biholomorphism from $A\ssmi A_1$ to $B\ssmi B_1$, with $A_1=q^{-1}(B_1)$.

As before there are three possible cases:  $f(\cbb)\cap B_1=\emptyset$,  $f^{-1}(B_1)$ is a discrete subset of $\cbb$ and $f(\cbb)\subset B_1$. If $f(\cbb)\cap B_1=\emptyset$, then we lift $f$ to the nonconstant 
descent datum $\psi(f)$ as in the Case 1 and prove that $\pis_0(\ycal,y)$ is not Brody hyperbolic. If $f(\cbb)\subset
 B_1$ then we repeat Case 3, but with all dimensions of the analytic spaces dropped, and recursively we will eventually end up with zero dimensional ramification, which is treated exactly as in the general case below, except being simpler.
Consider the case of  $f^{-1}(B_1)=:S_1$ being a discrete subset. 
Two things can happen: there exists an index $N$ such that for all $\nu>N$,   $f^n(\cbb)\subset B$ which is treated exactly as in the Case 2.
Otherwise, we proceed with few reduction steps. As usual we denote $F^\nu_{i(\nu)}$ a holomorphic map of the 
descent datum $\fsf^\nu$ on some disc, such that $\phi(\fsf^\nu)=f^\nu$. 
By the compactness of $\ycal$ we may assume $X$ had only a finite 
number of connected components (see beginning of Section \ref{aus}), hence, 
for a $\zeta\in S_1$, there is a punctured disc $\Delta_\zeta^*\subset \cbb\ssmi S_1$ and an index $k$ such that $F^j_{i(j)}(\Delta_\zeta^*)\subset X^{(k)}$ for infinitely many $j$. Considering only those indices $j$ we may assume 
all the images of $\Delta_\zeta^*$ through the sequence of functions $\{F^\nu_{i(\nu)}\}_n$ lie in the same connected component. 
Using again the compactness of $\ycal$, we may suppose that the sequence of points $\{F^\nu_{i(n)}(\zeta)\}_n$
converges to a point $a_1\in X^{(k)}$: otherwise it must necessarily converge to a point $a_1$ on the boundary of $X^{(k)}$ in $Z^{(k)}$ and $q^{-1}(q(a_1))$ has a point in the interior of another connected component $X^{(h)}$
of $X$, for which all the statements regarding $X^{(k)}$ hold as well. As far as our argument is concerned, we then replace $X^{(k)}$ with $X^{(h)}$. Let $U_{a_1}$ be an open neighborhood of $a_1$ such that the closure 
$\overline{U_{a_1}}\subset X^{(i)}$.  Then

{\em there is a closed simple path $\gamma$ in $\cbb\ssmi S_1$ around $\zeta$ and close enough to it for 
$$F^\nu_{i(\nu)}(\gamma([0,1]))\subset U_{a_1}$$ 
for all $\nu\geq N\gg 0$.}\\
Indeed, as in Section \ref{top}, let $H$ be a 
distance on $Q(\ycal)$ and $q^*H$ its pullback on $X$; set $d$ to be a positive real number such that the set of points in $X^{(i)}$ whose  ${\rm d}_{q^*H}$ distance from $a_1$ is less than $d$ is contained in $U_{a_1}$. Choose a closed simple path $\gamma$
as in the statement, in a way that ${\rm d}_H(f(\gamma([0,1])), f(\zeta))=d/2$. Then there exists $N$ such that  ${\rm d}_{q^*H}(F^n_{i(n)}(\gamma([0,1])), a_1)\leq d$ for all $n\geq N$, since $q(a_1)=f(\zeta)$ and $q$ is a contraction with respect to the length functions $q^*H$ and $H$.

Replacing $X$ with $\overline{U_{a_1}}$, we can consider the following simplified new situation: $X$ is a connected, compact  topological space, $Q(\ycal)$ compact, $A$, $A_1$, $B$ and $B_1$ all compact subspaces.

Give $X$ a structure of CW complex with finitely many cells such that $A$ and $A_1$ are sub CW complexes. Let
$l={\rm d}_H(f(\gamma([0,1])), B_1)$; by possibly refining the cell structure we may take the $q^*H$-sizes of the cells
to be of dimension ranging between $l/4$ and $l/2$. Call $V_A$ the sub CW complex $A\subset V_A\subset X$
made of cells $e_i$ having nonempty intersection with $A$.
\begin{enumerate}
\item[1)] Since the sizes ${\rm d}_{q^*H}(e_i)\geq l/4$ for all $i$, and $f(\gamma([0,1]))\subset B$, there exists an $N$
such that $F^n_{i(n)}(\gamma([0,1]))\subset V_A$ for all $n\geq N$;
\item[2)] because ${\rm d}_{q^*H}(e_i)\leq l/2$, there exists $M\geq N$ such that, if $e_j$ are the cells of $V_A$ with non
empty intersection with the image of $F^\nu_{i(\nu)}\circ \gamma$ for any of the $\nu\geq M$, then $e_j\cap A_1=\emptyset$.   
\end{enumerate} 
 
Finally, let $R\geq M$ with the property that ${\rm d}_{q^*H}(F^\nu_{i(\nu)}(\gamma([0,1])), A)<l/8$. Then we can remove appropriately chosen points $p_1,\cdots, p_k$ from the interiors of some cells in $V_A$ to obtain a topological
space $V_A'$ with a deformation retract $r: V_A'\to A$ such that $r(F^R_{i(R)}(\gamma([0,1])))\subset A\ssmi A_1$.

Call $\gamma'$ the closed path $r\circ F^R_{i(R)}\circ\gamma$. Then we have that 
$$
[f\circ \gamma]=q_*(F^R_{i(R)}\circ \gamma)=q_*(\gamma')\in \pi_1(A\ssmi A_1,a)
$$
thus $f$ lifts to $A$ on a small punctured disc $\Delta_\zeta^*$, being $[f\circ \gamma]$ in the image of 
$$
{q_{|A\ssmi\,A_1}}_*:\pi_1(A\ssmi A_1,a)\to \pi_1(B\ssmi B_1, q(a))
$$ 
and $q_{|A\ssmi\, A_1}$ being a topological covering. As in the Case 2, we lift 
$f_{|\cbb\ssmi\, S_1}$ to $X$ (actually to $A$) and then extent the local lifting $\widetilde{f}_\zeta$ defined on the puntctured disc $\Delta^*_\zeta$ to $\Delta_\zeta$. This argument produces a nonconstant descent datum in $\pis_0(\ycal,y)(\cbb)$. 

\medbreak

This concludes the proof of the part ``$\pis_0(\ycal,y)$ Brody hyperbolic implies $\pis_0(\ycal,y)$ Kobayashi hyperbolic'' part of the Theorem \ref{BRO}. It remains to be shown that, if $\pis_0(\ycal,y)$ is Brody hyperbolic, then
$\pis_1(\ycal,y)$ is Kobayashi hyperbolic. However, by the proof so far written, we see that the Brody hyperbolicity
implies that $c(\ycal,y)<+\infty$; by Lemma \ref{lcat1}, both of $\pis_0(\ycal,y)$ and $\pis_1(\ycal,y)$ are Kobayashi hyperbolic.
\enddemo
The proof of the Theorem \ref{BRO} has highlighted the connection between hyperbolicity of $Q(\ycal)$ as
a complex space and $\ycal$ as  Deligne-Mumford analytic stack.

\begin{cor}\label{equivalenze}
Let $X\to\ycal$ be a compact  Deligne-Mumford analytic stack. Then
\begin{enumerate}
\item[1)] if $Q(\ycal)$ is hyperbolic, $\ycal$ is hyperbolic;
\item[2)] $\ycal$ is hyperbolic if and only if the presheaf $\pis_0(\ycal,y)$ is hyperbolic if and only if ${\sf Ob}(\ccal
\xdot)$ is an hyperbolic presheaf.
\end{enumerate}
\end{cor}


\section{Hyperbolicity and the coarse moduli space}\label{hyp-coars}

While hyperbolicity of the coarse moduli space implies hyperbolicity of the stack, the converse is generally
not true. Let $p:X\to \ycal$ be an analytic stack with a one dimensional torus $\mathbb T$ as coarse moduli space $Q(\ycal)$
and such that the ramification of $q:X\to Q(\ycal)$ is at least $4$ in each of the points of the fiber over
two branch points $z_1$ and $z_2\in Q(\ycal)$. We show that such a stack $\ycal$ is hyperbolic by proving that no nonconstant map $f:\cbb\to Q(\ycal)$ can be lifted to $X$. Indeed, let $\tilde{f}:\cbb\to \cbb$ be the lifting of $f$ to the universal 
covering $h:\cbb\to Q(\ycal)$. Assume that $f$ admits local liftings to $X$. Since $f$ must be surjective, and $h$
is unramified, any two points $w_1$ and $w_2$ in the $h$ fiber of $z_1$ and $z_2$ respectively will have the following 
property: $\tilde{f}^{-1}(w_1)$ and $\tilde{f}^{-1}(w_2)$ are two sets of points on which $\tilde{f}$ have  ramification
divisible by $N_1$ and $N_2$, respectively, both of them greater than $4$. However, it is impossible for
such an entire function to exist. Indeed, otherwise there would be two functions $\alpha$ and $\beta$ such 
that $\alpha^{N_1}=\tilde{f}-w_1$ and $\beta^{N_2}=\tilde{f}-w_2$ and satisfying the equation 
$x^{N_1}-y^{N_2}=w_2-w_1$. The function $z\to (\alpha(z),\beta(z))$ is defined over $\cbb$ and  nonconstant
with image in the affine part of a curve of genus greater than $2$. This is impossible.  

\noindent The following is an example of such a stack: consider a covering of a torus $\mathbb T$ by a curve $q: X\to \mathbb T$
of degree $4$ and ramified in two points $x_1$, $x_2\in X$ with ramification index $4$, given by the composition 
$X\stackrel{q_2}{\to} X'\stackrel{q_1}{\to} \mathbb T$. Here both maps have degree $2$ and are ramified in only two points.
Their existence follows from
\begin{lemma}
Given a complex curve $C$, and $B=\{c_1, c_2,\cdots , c_{2k}\}\subset C$, there exists a degree $2$ 
covering $p:\widetilde C\to C$ with branch locus $B$.
\end{lemma}
\demo
It follows from the classification of the coverings of curves (see, for instance, \cite[ Prop. 4.9]{mir}). Coverings of $C$ 
with branch locus contained in $B$ are classified by homomorphisms $\pi_1(C\ssmi B, x)\to S_2$ with transitive image, where $S_2$ is the symmetric group on two elements. The condition on the ramification corresponds to 
the image of loops $\gamma_i$ around points of $B$ going to the transpositon in $S_2$. $\pi_1(C\ssmi B, x)$ can be presented as 
$$\langle \gamma_i, \alpha_j, \beta_j: \prod_i\gamma_i\cdot \prod_j
\alpha_j\beta_j\alpha_j^{-1}\beta_j^{-1}=1\rangle_{i=1,\cdots, 2k,\;\; j=1,\cdots, \mathrm{genus}(C)} 
$$
where $\alpha_j$ and $\beta_j$ are the generators of $\pi_1(C,x)$. Therefore, such an homomorphism exists
as $B$ contains an even number of points.
\hfill \enddemo

The covering $q$ is Galois, being the composition of two Galois extensions and since the monodromy group of
a covering map is isomorphic to the Galois group of the rational functions field extension of the curves involved (cfr. \cite{harris}, proposition at page 689, for instance). Therefore $[X/\mathrm{Gal}(\cbb(X)/\cbb(\mathbb T))]$ is a quotient stack with $\mathbb T$ as coarse moduli space satisfying the condition on the ramification stated above.


\subsection{Automorphisms of a compact  Deligne-Mumford analytic stack}

Stack hyperbolicity is expected to impose a peculiar behaviour to the stack itself, just like usual hyperbolicity does 
on a complex space. 
Here we prove that a compact  Deligne-Mumford analytic stack with a hyperbolic coarse moduli space  has few  automorphisms: 

 \begin{theorem}\label{auto}
Let $p:X\to \ycal$ be a compact Deligne-Mumford analytic stack with hyperbolic coarse moduli space. Then the isomorphism classes of the $2$-group
${\sf Aut}(\ycal)$ are finitely many. 
\end{theorem}

\demo By assumption the coarse moduli space $Q(\ycal)$ is hyperbolic and compact, hence it has a finite number of biholomorphisms (cf. \cite[Theorem 5.4.4]{KOB1}). We will prove that the kernel ${\sf Aut}_{Q(\ycal)}(\ycal)$ of the  $2$-homomorphism that associates to an automorphism $\ycal\to \ycal$ the induced biholomorphism $Q(\ycal)\to Q(\ycal)$ has a finite number of equivalence classes. Let $q:X\to Q(\ycal)$ be the canonical projection.

\begin{lemma}
There exist a dense open subspace $V\subset Q(\ycal)$ and
 a dense open substack $\ucal\subset\ycal$  
which is a gerbe over $V$. 
\end{lemma}
\demo An analytic stack $\wcal$ is a gerbe on a substack where the inertia stack morphism $\Delta':
\ical\to \ycal$ is flat. Consider the $2$ commutative diagram 
$$
\xymatrix{Z\subset W\ar[r]\ar[d]^f & \ical\ar[r]\ar[d]^{\Delta'} & \ycal\ar[d]^\Delta\\
f(Z)\subset X\ar[r] & \ycal \ar[r]  & \ycal\times\ycal}
$$
where the squares are cartesian. By a theorem of Frisch (cfr. \cite{fri}), $f$ is flat outside an analytic subset $Z\subset W$. As
$f$ is proper ($\Delta$ is by assumption, see \ref{procpt}),  $f(Z)$ is analytic in $X$. We conclude that $f:\ical
\times_\ycal X\ssmi f(Z)\to X\ssmi f(Z)$ is flat. Now, $X\ssmi f(Z)$ is $X_1=X\times_\ycal X$ invariant, as it is maximal
in the codomain onto which  $f$ is flat. Therefore, there exists an associated open dense substack $\ucal\subset 
\ycal$ with the property that the basechanged morphism $\ical_\ucal\to \ucal$ is flat, by faithfully flat descent.
This shows that $\ucal$ is a gerbe onto $V:=Q(\ycal)\ssmi q(f(Z))$. 
\enddemo

We also have that $\pi_1(V)$ is finitely generated, since the following holds

\begin{prop}
Let $C$ be a compact complex space and $B\subset C$ analytic; then  $\pi_1(C\ssmi B)$
is finitely generated. 
\end{prop}
\demo
Embedding $C$ in some $\R^N$ (cfr. \cite{toto}) and applying the triangularization 
theorem of \L ojasiewicz \cite[Theorem 3]{lo} we are reduced to the case where $C$ is a finite simplicial 
complex $K$ and $B$ is a subcomplex $H< K$.  Let $\gamma:S^1\to |K|\ssmi |H|$ be a closed path; we will show that $\gamma$ is homotopic to the topological realization of a simplicial map $\hat{\gamma}:\hat{S}^1\to |D|\subset
|K|\ssmi |H|$ for some simplicial complex $D<K$. This finishes the proof, since $D$ is a finite subcomplex, thus 
$\pi_1(D)$ is finitely generated. 
Let $\gamma$ be given. If the image of $\gamma$ does not intersect the star of $H$, then we can let 
$\hat{\gamma}$ to be a simplicial approximation of $\gamma$ in $K$. If $\gamma(S^1)$ intersects non 
trivially the star of $H$, it is possible to continously deform $\gamma$ inside each simplex of $\overline{|St(H)|}$
to get a $\gamma'$ homotopic to $\gamma$ with image not intersecting the star, but rather the link of $H$ in $K$.
In particular, the image of $\gamma'$ will be contained in the simplicial complex $D=K\ssmi St(H)$. Take now 
a simplicial approximation of $\gamma'$ in $D$.
\enddemo

We show that the restriction functor ${\sf End}_{Q(\ycal)}(\ycal)\to {\sf End}_V(\ucal)$ is injective on the associated 
isomorphism classes and then that the isomorphism classes of the codomain are finite.
Assume that two morphisms $f,g:\ycal\to \ycal$ are isomorphic over $\ucal$ and let $h:\ical\to \ycal$ be the pullback of  the diagonal morphism $\ycal\to \ycal\times_{Q(\ycal)}\ycal$ through
$(f,g):\ycal\to \ycal\times_{Q(\ycal)} \ycal$. 
Consider the diagram:
$$
\xymatrix{W\ar[r]\ar[d]^\phi & \ical\ar[rr]\ar[d]^{h} & & \ycal\ar[d]^\Delta\\
X\ar[r] & \ycal \ar[rr]^{(f,g)\hskip 25pt}  & & \ycal\times_{Q(\ycal)}\ycal}
$$
where the vertical morphisms are base changements of the diagonal $\Delta$.
Sections to $h$ correspond to 
isomorphisms between $f$ and $g$. The isomorphism between $f$ and $g$ over $\ucal$ yields a section 
$s:\ucal\to \ical$ to $h$. Because of the compactness of $\ycal$, the diagonal morphism is proper (see Proposition \ref{cptdiag}), $h$ is proper. Moreover, $\ycal$ being Deligne-Mumford, we conclude that its diagonal is \'etale, so $h$  is as well.  The same properties are inherited by the base changed 
holomorphic map $\phi:\ical\times_\ycal X=:W\to X$ of $h$. 
Let $a\in A=X\ssmi U$ be a point and $a\in V_a\subset X$ be a simply connected open neighborhood. By assumption
we have a section $s$ to $\phi$ over $U\subset X$. As $\phi$ is a topological covering, $s$ extends over $V_a$ 
because it is simply connected and of the unicity of the lifting, once we fixed base points. Thus, $s$ extends to
a section of $\phi$ over the whole of $X$, so that we get a map $s':X\to \ical$ making the diagram commute. $s'$
induces a stack morphism $[X_1\rightrightarrows X]\to \ical$ corresponding to a section $\ycal\to \ical$ through
the canonical equivalence $p: [X_1\rightrightarrows X]\to \ycal$.

We now prove that ${\sf End}_{V}(\ucal)$ has a finite number of isomorphism classes. ${\sf End}_{V}(\ucal)$ is a finite
disjoint union of gerbes: locally $\ucal$ is $V_i\times BG_i$ for some finite group $G_i$, and ${\sf Hom}(BG_i,BG_i)=
[{\sf Hom}(G_i,G_i)/G_i]$, the quotient stack of the finite set ${\sf Hom}(G_i,G_i)$ by conjugation (cfr. 
\cite[Proposition 5.3.5]{giraud}), therefore a finite disjoint union of $B{\rm Stab}_{i_j}$, where Stab$_{i_j}$ is the stabilizer of the $j$-th orbit of the $G_i$ action on ${\sf Hom}(G_i,G_i)$. Each such gerbe has a global object over $V$, namely the identity morphism ${\rm id}_\ucal$ of $\ucal$. It follows that 
${\sf End}_V(\ucal)$ is equivalent to $B_VG$ where $G={\sf Aut}({\rm id}_\ucal)$. $G$ is a Lie group endowed by a holomorphic map
$G\to V$ whose fibers are finite groups. Therefore, it suffices to show that 
there are a finite number 
of $G_V$-torsors on $V$ modulo isomorphism. Because of the finiteness condition on  $G\to V$, holomorphic 
$G_V$-torsors are the same as topological $G$-torsors, therefore $B_VG(V)/\text{iso}={\sf Hom}(\pi_1(V), G)$
is finite since $\pi_1(V)$ is finite.
 
\hfill \enddemo



\end{document}